\documentclass[12pt,leqno,fleqn,epsfig]{article}
\usepackage{amssymb,epsfig, amsmath, amsthm}
\usepackage{mathrsfs}       

\newcommand{\R}{\mathbb{R}}

\newcommand{\N}{\mathbb{N}}

\newcommand{\bbf}{\boldsymbol f}
\newcommand{\bg}{\boldsymbol g}

\newcommand{\bu}{\boldsymbol u}
\newcommand{\bv}{\boldsymbol v}
\newcommand{\bw}{\boldsymbol w}

\newcommand{\bA}{\boldsymbol A}
\newcommand{\bB}{\boldsymbol B}

\newcommand{\bD}{\boldsymbol D}
\newcommand{\bE}{\boldsymbol E}
\newcommand{\bF}{\boldsymbol F}

\newcommand{\bI}{\boldsymbol I}
\newcommand{\bJ}{\boldsymbol J}

\newcommand{\bL}{\boldsymbol L}

\newcommand{\bP}{\boldsymbol P}

\newcommand{\bT}{\boldsymbol T}

\newcommand{\bW}{\boldsymbol W}

\newcommand{\bfvarphi}{\boldsymbol\varphi}

\newcommand{\bfpsi}{\boldsymbol\psi}
\newcommand{\bfphi}{\boldsymbol\phi}

\newcommand{\bfzeta}{\boldsymbol \zeta}

\newcommand{\bPhi}{\boldsymbol \Phi}

\newcommand{\mes}{\operatorname{\rm meas}}    
\newcommand{\esssup}{\operatorname*{ess\,sup}}

\newcommand{\const}{\operatorname*{const}}

\newcommand{\be}{\begin{equation}}
\newcommand{\ee}{\end{equation}}
\newcommand{\bea}{\begin{eqnarray}}
\newcommand{\eea}{\end{eqnarray}}
\newcommand{\bean}{\begin{eqnarray*}}
\newcommand{\eean}{\end{eqnarray*}}

\newcommand{\intl}{\int\limits}
\newcommand{\suml}{\sum\limits}

\newcommand{\Beweisende}{\rule{0.2cm}{0.2cm}}

\newcommand{\D}{\displaystyle}

\newcounter{secnum}

\newtheorem{thm}{Theorem}[section]
\newtheorem{cor}[thm]{Corollary}
\newtheorem{lem}[thm]{Lemma}

\theoremstyle{definition}

\newtheorem{defin}[thm]{Definition}
\newtheorem{rem}[thm]{Remark}

\title{On the local pressure of the Navier-Stokes equations and related systems}
 \author{J\"{o}rg Wolf \\
\ \\
 Department of Mathematics\\
Humboldt University of Berlin\\
Unter den Linden 6, 10099 Berlin, Germany\\
 e-mail: jwolf@math.hu-berlin.de
 \\
\date{}}

\begin{document}
\maketitle

\begin{abstract}
In the study of local regularity of weak solutions to systems related to incompressible viscous fluids local energy estimates serve as important ingredients. 
However, this requires certain  informations on the pressure. This fact has been used  by V. Scheffer in the  notion of a suitable weak 
to the Navier-Stokes equation, and  in the proof of the partial regularity due to   Caffarelli. Kohn and Nirenberg.  
In general domains, or in case of  complex viscous fluid models a global pressure doesn't necessarily exist. 
To overcome this problem, in the present paper we construct a local pressure distribution by showing  that every 
distribution $ \partial _t \bu +\bF $, which vanishs on the set of smooth solenoidal vector fields can be represented by  a distribution  $ \partial _t \nabla p_h +\nabla p_0 $, where $\nabla p_h \sim \bu  $ and $ \nabla p_0 \sim \bF$. 
\end{abstract}

\vspace{0.5cm}  

\vspace{0.5cm}  

\tableofcontents

%%% ----------------------------------------------------------------------
%       SECTION 1
%%% ----------------------------------------------------------------------
\section{Introduction}
\label{sec:-1}
\setcounter{secnum}{\value{section} \setcounter{equation}{0}
\renewcommand{\theequation}{\mbox{\arabic{secnum}.\arabic{equation}}}}

Let $\Omega \subset \R^n$ $(n\ge 2)$ be an open set.  By  $W^{ 1,\, q}_0(\Omega )^n$ 
$(1<q<+\infty )$  we denote the closure of $C^\infty _{\rm c} (\Omega )^n$  under the usual Sobolev 
norm. By $W^{ - 1,\, q}(\Omega )^n$  we denote the dual of $W^{ 1,\, q'}_0(\Omega )^n$ \,\footnotemark\,.  
For $p\in  L^q(\Omega )$  by  $\nabla _q p $ we mean the functional in 
$W^{ -1,\, q}(\Omega )^n$ 
determined by 
\footnotetext{\,Here $q' = \frac {q} {q-1}$ if $1<q<+\infty $, $q'=+\infty $ if $q=1$ and $q'=1$ if $q=+\infty $.}
\be\label{1.1}
\langle \nabla _q p, \bfpsi  \rangle= -\intl_{\Omega} p \nabla \cdot \bfpsi dx,\quad 
\bfpsi\in  W^{ 1,\, q'}_0(\Omega )^n \,\footnotemark\,.   
\ee
\footnotetext{\,Here $\nabla \cdot \bv$ stands for the divergence 
$ \suml_{i=1}^{n} \partial_i v^i$ of the field $\bv= (v^1, \ldots, v^n)$.}
In the present paper we are interested in the existence of a projection 
$\bE_q:  W^{ 1,\, q'}_0(\Omega )^n \rightarrow W^{ 1,\, q'}_0(\Omega )^n $ such that 
\begin{align}
 \hspace*{0.5cm}  \bE_{q'}  (\bfpsi ) &={\bf 0} \quad &&\forall\, \bfpsi \in C^\infty _{\rm c}(\Omega )^n
\,\, \mbox{with}\,\, \nabla \cdot \bfpsi =0,
\hspace*{5cm}  
\label{1.2}
\\
\bE_{q'}  (\nabla\phi ) &= \nabla \phi  \quad &&\forall\, \phi \in C^\infty _{\rm c}  (\Omega ).  
\label{1.3}
\end{align}
As we will see below, the existence of such projection is ensured if $\Omega $  is sufficiently regular. 
However,  \eqref{1.2} and \eqref{1.3} do not guarantee the uniqueness of $ \bE_{ q'} $, so we may replace condition \eqref{1.3}
by a more restrictive one.  Moreover, due to \eqref{1.2} the dual projection  
$\bE_q^\ast: \bW^{ -1,\, q}(\Omega )^n \rightarrow 
\bW^{ -1,\, q}(\Omega )^n $ enjoys the property 
\be\label{1.4}
\hspace*{0.5cm}  \bE^\ast_q(\nabla _q p) =\nabla _q p\quad \forall\, p\in L^q(\Omega ),
\ee
 so that $\bE^\ast_q $  appears to be a useful tool for constructing the pressure of weak solutions to the equations 
modelling the motion of an incompressible fluid.  In three dimensions,  this systems consist of four equations, formed by the conservation of momentum and the conservation 
of volume including four unknowns, the pressure $p$ and the velocity field $\bu=(u^1, u^2, u^3)$.  
Among this fluid models perhaps  the  Navier-Stokes system is one of the most popular.  Since  the pioneering work by J.\,Leray \cite{Leray1934}  the theory of the Navier-Stokes equations has been 
widely developed, where fundamental problems such as existence of weak solutions, conditions for global and local  regularity or asymptotical  behaviour have been solved. However, despite strong efforts one of the most important question, the existence of a unique  global regular solution for general smooth data is still open. Partial answers to this fundamental  question have been given, such as sufficient conditions for global or local regularity,  which have been relaxed step by step in recent years.  Concerning the local regularity, a first result goes back to V.\,Scheffer \cite{Scheffer1976}, who introduced the notion of a  
{\it suitable weak solutions} to the Navier-Stokes equations,  fulfilling  a local energy inequality. 
In 1982,  based on Scheffer's notion, L.\,Caffarelli, R.\,Kohn and L.\,Nirenberg  \cite{CaKoNi1982} obtained 
an optimal result of partial regularity of suitable weak solutions,  by showing that the one-dimensional parabolic Hausdorff measure of the singular set is zero.  Recently, this result has been improved logarithmically  by 
Choe and Lewis in \cite{ChLe2000}.   For  alternative  proofs of  Caffarelli-Kohn-Nirenberg theorem we refer  to  \cite{Lin1998, LaSe1999, Vasseur2008, Wolf2008}.    

\vspace{0.2cm}
 \hspace*{0.5cm}
The proof of the Caffarelli-Kohn-Nirenberg theorem 
rests  on decay estimates, derived  from the local energy inequality, which holds for 
suitable weak solutions (cf. \cite{Scheffer1976, CaKoNi1982}).   
Unlike weak Leray-Hopf solutions, which are constructed by using Galerkin approximation,  
suitable weak solutions have to be  constructed differently. 
It is still unclear, whether a Leray-Hopf 
solution is suitable or not. Furthermore, since the existence of a suitable weak 
solution depends on the existence of a global pressure (cf. \cite{SoWa1986}, \cite{Sohr2001} and \cite{FaKoSo2005} for general uniform $C^2$-domains)  this method requires that $\Omega $ is sufficiently regular.  Note, that the same problem occurs in other fluid models, such as non-Newtonian fluids and fluids with variable viscosity. In the recent paper  \cite{Wolf2014a} a local  pressure projection has been introduced to obtain the 
Caffarelli-Kohn-Nirenberg theorem for local suitable weak solutions to the Navier-Stokes equations in  arbitrary cylindrical domain domains $ Q= \Omega \times [0,T]$, for any open set   $ \Omega \subset \R^{3}$.  
In fact, this method suggests to even work  with distributional weak solutions to the 
Navier-Stokes equations with variable viscosity, and related systems. However, one has to be careful as the following example of a potential-like solution  shows (cf. \cite[footnote p.79]{GaRio1983}). 

\vspace{0.3cm}
{\bf Example}. Let $\eta : ]0,T[ \rightarrow \R$ be any function and let $\phi : \Omega \rightarrow \R^3 $ be harmonic.  Then it is not difficult to check that $\bu(x,t)= \nabla \phi (x) \eta (t) $ solves 
the Navier-Stokes equations
\begin{align*}
 \nabla \cdot \bu &=0 \quad\mbox{in}\quad \Omega \times ]0,T[, 
\\     
\partial _t \bu + (\bu\cdot \nabla) \bu - \Delta \bu&= - \nabla p
\quad\mbox{in}\quad \Omega \times ]0,T[
\end{align*}
in the sense of distribution, where the pressure is given by the following distribution
\[
p= -\phi \eta' - \frac {1} {2} |\nabla \phi |^2 \eta^2.    
\]
Note that in this example we have not imposed any boundary condition on $ \bu $. In fact, in case of no slop boundary 
condition and if $ \Omega $ is sufficiently smooth $ \bu $ becomes trivial.

\hspace{0.5cm}
As the above example shows, the pressure  might not be a  Lebesgue function and it is unlikely 
to improve the time regularity for a distributional solution $\bu$.   On the other hand,   
the pressure has the following form  
\be\label{i0}
p=\partial _t p_h + p_0,\quad p_h =- \phi \eta,\quad p_0 = -\frac {1} {2} |\nabla \phi |^2 \eta^2, 
\ee
where $p_h(t)$ is harmonic for a.\,e. $t\in  ]0,T[$,  which suggests to introduce the local 
pressure taken as in  \eqref{i0} on suitable subdomains.  This will be done  with help of a  projection 
$\bE_q^\ast$ fulfilling \eqref{1.2} and \eqref{1.3}.   In fact,  such  method of  pressure representation on subdomains 
has introduced first in \cite{Wolf2007a}  and  later used in \cite{DiRuWo2010}  to construct a weak solution to the equations of non-Newtonian fluids in general domains.  
This method  has played also an significant role to  achieve  further  results  concerning 
existence and regularity of weak solutions to  models related to incompressible viscous fluids 
(cf. \cite{BuFeNeWo2008},  \cite{Wolf2008}, \cite{Wolf2010}).

\vspace{0.2cm}
\hspace*{0.5cm} 
In  the present paper we wish to generalize the method introduced above such that it can be used  for any given distributional solution to the Navier-Stokes equation 
or related systems in a cylindrical domain $ Q= \Omega \times ]0,T[$ \,$ (0<t<+\infty)$ governing incompressible viscous fluids. Our main result will be the characterization of distributions of the form   $\partial _t \bu + \nabla \cdot \bA$  in  $ Q$ vanishing on the space of all smooth solenoidal fields 
with compact support in $Q$, by a distribution involving gradient fields only. More precisely,   
 for every $C^1$ subdomain $G\Subset \Omega $ \,\footnotemark\,  there are pressure functions $p_{h, G}(t)$ and $p_{0,G} (t)$ with  $\nabla p_{ h,G} (t)\sim \bu(t)$  and $ p_{ 0,G}(t) \sim \bA(t) $  (for a.\,e. $t\in  ]0,T[$) 
satisfying  $\partial _t \bu + \nabla \cdot \bA= -\partial _t \nabla p_{h, G}- \nabla p_{0, G} $ 
in  the sense   of distributions, i.\,e. 
\footnotetext{\,Here, for two sets $A, B \subset \R^n$,  the notation $A\Subset B$ 
means $\overline{A }\subset B$, and $ \overline{A} $ compact. }
 \be\label{i13}
\intl_{0}^{T}\intl_{G}  -\bu \cdot \partial _t\bfvarphi  - \bA: \nabla \bfvarphi  dx   d t 
= \intl_{0}^{T}\intl_{G}  - p_{h,G}  \cdot \partial _t \nabla \cdot \bfvarphi + 
p_{0,G}  \nabla \cdot \bfvarphi   dx   d t \,\footnotemark\,
\ee 
for all  $\bfvarphi \in  C^\infty _{\rm c}(Q) ^n$.
\footnotetext{\, For matrices $\bA, \bB\in  \R^{n\times n}$ by $\bA:\bB$ we denote the scalar product 
$\suml_{i, j=1}^{n}A_{ij}B_{ij} $.  } Here $\bu = (u^1, \ldots , u^n)$ stands for the velocity field of the fluid and $\bA= \{A_{ij} \}$ for an  $n\times  n$ tensor modelling the fluid system.   For instance, 
 the Navier-Stokes equation  is   modelled by 
\[
\bA = \bu\otimes \bu - \nu \nabla \bu,\quad \nabla \cdot \bu =0 \quad \mbox{in}\quad Q,
\]
where $\nu =\const >0$ denotes the viscosity of the fluid.  

\vspace{0.2cm}
\hspace*{0.5cm} 
The paper is organized as follows. In Section\,2 we provide some notations and function spaces 
used throughout the paper.  In Section\,3 we introduce the space $G^{-1,\,q}(\Omega )^n $ containing all functionals $\bu^\ast \in  W^{ -1,\, q}(\Omega )^n$ which vanish on the space of solenoidal vector fields.  As we will see below the space $ {G}^{-1,\,q}(\Omega )^n $  contains distributions of the form $\nabla p$.    Then we are interested in domains $\Omega $ for which there exists  a projection $\bE_q^\ast$  from 
$W^{ -1,\, q}(\Omega )^n$ onto $ G^{-1,\,q}(\Omega )^n $  being the dual of 
 a projection $\bE_{q'} $ of  $W^{ 1,\, q'}_0(\Omega )^n$  onto a closed subspace, fulfilling 
\eqref{1.2} and \eqref{1.3}.  Such domains will be called $\nabla_q $-regular, 
and they will be used for the construction of the pressure 
representation $\nabla p = \partial _t \nabla p_{h}+  \nabla p_0$.  
Section\,4 deals with the existence of unique $q$-weak solutions to the Stokes-like system, 
which  appears to be a  sufficient criterion for $\nabla _q$-regularity. 
This will be verified for the 
following cases,  (i) $\Omega =\R^n$, (ii) $\Omega $ is a bounded $C^1$ domain and (iii) 
$\Omega $ is a exterior $C^1$-domain.  
In Section\,5 we  first introduce a sufficient criterion for $\nabla _q$-regularity based on the 
existence and uniqueness of $q$-weak solutions to the Stokes-like system. 
Next, in Section\,6 we present our first main result concerning the local pressure decomposition 
for time-dependent distribution.  Then, in  Section\,7 for the case $ q=2$ we provide a global pressure representation for a general 
domain by using orthogonal projections due to the Hibert space structure. We complete our discussion by applying the result of Section\,7 to the case of the generalized Navier-Stokes equations.    
At the end of the paper we have added  Appendix\,A, recalling some well-known properties of vector-valued function in suitable Bochner spaces and Appendix\,B, where we discuss the continuity of time dependent potentials of spatial gradients being continuous in space and time, 
which is related to the continuity of the harmonic pressure.

%%% ----------------------------------------------------------------------
%       SECTION 2
%%% ----------------------------------------------------------------------
\section{Notations and function spaces}
\label{sec:-2}
\setcounter{secnum}{\value{section} \setcounter{equation}{0}
\renewcommand{\theequation}{\mbox{\arabic{secnum}.\arabic{equation}}}}

Let $\Omega \subset \R^n$ $(n \in \N,  n\ge 2)$ 
denote a domain.  If necessary the properties of $\Omega $ will be specified. By $W^{k,\, q}(\Omega ), 
W^{k,\, q}_0(\Omega )$ \, $(1\le q < +\infty; k \in  \N)$  
we denote the usual Sobolev spaces.  If $\Omega $ is bounded, the space  $ W^{k,\, q}_0(\Omega )$ will 
be equipped with the norm 
\[
\|u\|_{ W^{k,\, q}_0} = \Big(\suml_{|\alpha |=k} \|D^\alpha u\|^q_{L^q}\Big)^{1/q}, 
\]
otherwise with the usual Sobolev norm. 
For $1<q< +\infty$,   the dual of  $ W^{k,\, q'}_0(G)$ \,\footnotemark\, 
will be denoted by $W^{ -k,\, q}(\Omega )$.
Throughout, without any reference vector valued or tensor-valued functions will be denoted  by boldface letters. 
%\footnotetext{\,Here $q' = \frac {q} {q-1}$ if $1<q<+\infty $, $q'=+\infty $ if $q=1$ and $q'=1$ if $q=+\infty $.}

\vspace{0.2cm}
\hspace*{0.5cm} 
Next,  by $C^\infty _{\rm c, div} (\Omega )^n$ we denote the space of all smooth solenoidal vector fields 
$\bfpsi : \Omega \rightarrow \R^n$ having its support in $\Omega $.  Then by 
$W^{ 1,\, q}_{0, \rm div} (\Omega )^n$ we denote the closure of  
$C^\infty _{\rm c, div} (\Omega )^n$  with respect to the norm in $W^{ 1,\, q}_0(\Omega )^n$
  \, $(1\le q < +\infty )$.  Similarly, by  $L_{\rm div}^q (\Omega )^n$ we denote the closure of  
$C^\infty _{\rm c, div} (\Omega )^n$  with respect to the $L^q$-norm  \, 
$(1\le q< +\infty )$. 
Furthermore, in case $\mes \Omega <+\infty $ by  $L^q_0(\Omega )$ we denote the subspace 
of all $p\in  L^q(\Omega )$ such that 
$\intl_{\Omega }  pdx=0 $. 

\vspace{0.1cm}
\hspace*{0.5cm} 
Let $\bJ_q=\bJ_{q, \Omega }: W^{ 1,\, q}_0(\Omega )^n \rightarrow   W^{ -1,\, q}(\Omega )^n $ 
be defined by  
\[
\begin{cases}
{\D \langle \bJ_q \bu, \bv\rangle = \intl_{\Omega } \nabla  \bu: \nabla \bv  dx  }  \,\, & \hspace*{-1cm}\mbox{if $\Omega $ is bounded}\quad 
\\[0.2cm]
{\D \langle \bJ_q \bu, \bv\rangle = \intl_{\Omega } \nabla  \bu: \nabla \bv + \bu\cdot \bv dx  }  \,\, &\hspace*{-1cm}\mbox{if $\Omega $ is unbounded}\quad
\\[0.2cm]
\bu \in  W^{ 1,\, q}_0(\Omega )^n, \bv \in  W^{ 1,\, q'}_0(\Omega )^n\quad (1<q<+\infty ).
\end{cases}
\]
Note that $\bJ_{q}$ defines an isomorphism in the cases (i) $\Omega =\R^n$, (ii)  
$\Omega =\R^n_+$ or  (iii) $ \Omega $ is a  $C^1$-domain 
with compact boundary. Note that if $ \Omega \subset \R^{n}$   
is a bounded Lipschitz domain then there exists $ 3<q_1<+\infty$  if $ n \ge  3$  
or $4< q_0<+\infty $   if $ n =2$ such that $ \bJ_{q}$ is an isomorphism for all 
$ q_0'<q<q_0$, which is in some sense sharp (cf.  \cite{JeKe1989}).  
In the special case  $q=2$, $\bJ_{2}$  defines isomorphism for any domain,  
which is due to the Hilbert space structure.   
In fact,  here  $\bJ_{2} $ coincides  with the duality map for $W^{ 1,\, 2}_0(\Omega )^n$, 
while $\bJ_{2 } ^{-1} \bu^\ast\in  W^{ 1,\, 2}_0(\Omega )^n$  appears to be the 
Riesz representation of the functional $\bu^\ast \in W^{ -1,\, 2}(\Omega )^n $.  

\vspace{0.1cm}
\hspace*{0.5cm} 
Let $\bu \in  L^1_{\rm loc}(\Omega )^n$. 
We say $\bu \in  W^{ -1,\, q}(\Omega )^n$  for $1<q<+\infty $  if there exists 
$c=\const>0$ such that 
\[
\langle \bu, \bfpsi  \rangle  := 
\intl_{\Omega}  \bu \cdot \bfpsi  dx  \le c \|\bfpsi \|_{W^{1,q'} } \quad \forall\, \bfpsi\in C^\infty _{\rm c}(\Omega )^n.  
\]
Hence, there exists a unique $ \bbf \in W^{ -1,\, q}(\Omega )^n$ such that 
$ \langle \bu, \bfpsi  \rangle = \langle \bbf, \bfpsi  \rangle$  for all $ \bfpsi \in C^\infty _{\rm c}(\Omega )^n$. In this case we may identify  $ \bu $ with  $ \bbf $, which justifies the above notation.  

\vspace{0.1cm}
\hspace*{0.5cm}  
For our discussion below we use the notation $\nabla _q = \nabla_{q, \Omega } $  for the gradient operator, mapping from $L^q(\Omega )$ into $W^{ -1,\, q}(\Omega )^n$ defined by \eqref{1.1}.  
From this definition we immediately derive  that the dual $\nabla ^\ast_q $  equals the divergence operator 
$- \nabla \cdot$  mapping from $W^{ 1,\, q'}_0(\Omega )^n$ into $L^{q'} (\Omega )$, i.\,e. 
$\nabla ^\ast_q \bv = -\nabla \cdot \bv$ for $\bv \in W^{ 1,\, q'}_0(\Omega )^n$. 

%As will see in the next section in general we can't expect that the range of $\nabla _q$ is closed, 
%especially,  if $\Omega $ is unbounded.  This motivates to define the space $G^{-1,\,q}(\Omega )^n $ 
%containing all $\bu^\ast\in  W^{ -1,\, q}(\Omega )^n$ represented by a function $f \in  
%L_{\rm loc}^1(\Omega )$ which will be the subject of discussion of the next section. 

\vspace{0.2cm}
\hspace*{0.5cm} 
Let  $X$ be a Banach space with norm  $\|\cdot \|_X$.  
Let $-\infty \le a < b \le +\infty $. 
By $L^s(a,b; X)$ ($1\le s \le +\infty$) we denote the   space of 
all Bochner measurable functions  $f: ]a,b[ \rightarrow X$ such that 
\[
\intl_{a}^{b} \|f(t)\|^s _X d t < \infty \,\,\, \mbox{if} \,\,1\le s< \infty; \quad
\esssup_{t\in (a,b)}  \|f(t)\|_{X} < \infty \,\,\,\mbox{if}\,\,q=\infty.   
\]

%%% ----------------------------------------------------------------------
%       SECTION 3
%%% ----------------------------------------------------------------------
\section{The space $G^{-1,\,q}(\Omega )^n$}
\label{sec:-3}
\setcounter{secnum}{\value{section} \setcounter{equation}{0}
\renewcommand{\theequation}{\mbox{\arabic{secnum}.\arabic{equation}}}}

The present section deals with  properties of functionals  $\bu^\ast\in  W^{ -1,\, q}(\Omega )^n$  
vanishing  on $C^\infty _{\rm c, div}(\Omega )^n $.   To this end,  we introduce  the notion 
of $W^{ -1,\, q}$-potential. 

\begin{defin}
\label{def3.1}
Let $1<q<+\infty $.  
A function $p \in  L^q_{\rm loc}(\Omega ) $ is called a $W^{ -1,\, q}$-{\it potential} if there exists
$\bu^\ast\in  W^{ -1,\, q}(\Omega )^n$ such that 
\be\label{3.1}
\langle \bu^\ast, \bfpsi  \rangle = - \intl_{\Omega} p \nabla \cdot 
\bfpsi dx  \quad \forall\, \bfpsi\in C^\infty _{\rm c}(\Omega )^n.
\ee 
The set of all $W^{ -1,\, q}$-potentials will be denoted by $L^q_{\rm pot}(\Omega )$.  

\hspace*{0.5cm} 
If  $p\in L^q_{\rm pot}(\Omega )$ and $\bu^\ast\in  W^{ -1,\, q}(\Omega )^n $ fulfilling 
\eqref{3.1} we use the brief notation $\bu^\ast = \nabla _q p$.  
Then we define 
\begin{equation}
G^{-1,\,q}(\Omega )^n := \{\nabla _q p\in  W^{ -1,\, q}(\Omega )^n  \,|\, p\in  L^q_{\rm pot}(\Omega ) \}.
\label{3.1a}
\end{equation}
\end{defin}

\begin{rem}
\label{rem3.2}
If   $p\in  L^q(\Omega )$ we have $-\intl_{\Omega} p \nabla \cdot \bfpsi dx   \le n\|p\|_{L^q} 
\|\bfpsi \|_{W^{ 1,\, q'}}$ for all $\bfpsi \in  C^\infty _{\rm c}(\Omega )^n $,  which shows that $p\in  L^q_{\rm pot} (\Omega )$.  Thus, there exists a unique 
$\bu^\ast\in  W^{ -1,\, q}(\Omega )^n$ such that $\bu^\ast = \nabla _q p$. 
This implies  
\be\label{3.2}
L^q(\Omega )\subset  L^{q} _{\rm pot}(\Omega)\subset L^q_{\rm loc}(\Omega ). 
\ee
\end{rem}

\hspace*{0.5cm}  
The next lemma provides  a well-known characterization of  $G^{-1,\,q}(\Omega )^n $ 
(see also \cite[Cor.III.5.2]{Galdi1994} ).

\begin{lem} 
\label{lem3.3}
For $\bu^\ast\in  W^{ -1,\, q}(\Omega )^n$ the following statements are equivalent 

\vspace{0.3cm}
$1^\circ$ $\langle \bu^\ast, \bfpsi \rangle =0$ for all  $\bfpsi \in  
W^{ 1,\, q'}_{0, \rm div} (\Omega )^n$;

$2^\circ $  $\, \bu^\ast \in  G^{-1,\,q} (\Omega )^n$.
\end{lem}

{\bf Proof}:  1.  The implication $2^\circ  \Rightarrow 1^\circ $ holds, since $ C^\infty _{\rm c, div} (\Omega )^n $
is dense in $W^{ 1,\, q'}_{0, \rm div} (\Omega )^n$.  

\hspace*{0.5cm}  
2.  To prove $1^\circ  \Rightarrow 2^\circ$ we choose a sequence of balls  $B_i \in \Omega$ $(i\in \N)$ 
such that $\Omega_j := \cup_{i=1}^j B_i $  is connected for all $j\in  \N$, and $\cup_{j=1} ^\infty \Omega _j = \Omega $. Since $W^{ 1,\, q'}_{0, \rm div} (\Omega _j)^n \hookrightarrow W^{ 1,\, q'}_{0, \rm div} (\Omega)^n $ there holds 
\[
\langle \bu^\ast, \bv \rangle  =0\quad \forall\, \bv\in W^{ 1,\, q'}_{0, \rm div} (\Omega _j)^n. 
\]
By the aid of \cite[Cor.III.5.1]{Galdi1994} we get  a unique $p_j \in L^q(\Omega _j)$ with $(p_j)_{B_1} =0$ and 
\[
\langle \bu^\ast, \bv \rangle = \intl_{\Omega_j} p_j \nabla  \cdot \bv dx\quad \forall\, 
\bv \in W^{ 1,\, q'}_0(\Omega _j)^n\quad (j\in  \N).   
\]
As $ W^{ 1,\, q'}_0(\Omega _j)^n \hookrightarrow W^{ 1,\, q'}_0(\Omega _{j+1} )^n$  we see that 
$p_{j+1} - p_{j} = \const$ a.\,e. in $\Omega _j$, which must vanish since $( p_{j+1} - p_{j})_{B_1} =0$. 
Thus, $p_j = p_{j+1}|_{\Omega _j} $.  Hence, there exists a unique  $p\in L^q_{\rm loc}(\Omega ) $  
such that  $p_{B_1}=0$ and $p|_{\Omega _j}= p_j $.  Thus, setting $p=p_j$ a.\,e. in $\Omega _j$ $(j \in \N)$  
it follows  $p\in  L^q_{\rm pot}(\Omega ) $  satisfying 
$\nabla _q p =\bu^\ast$. Whence, $ \bu^\ast \in  G^{-1,\,q} (\Omega )^n$. 
\hfill \Beweisende

\begin{rem}
\label{rem3.4}
 Lemma\,\ref{lem3.3} implies that $G^{-1,\,q}(\Omega )^n = 
\Big(W^{ 1,\, q'}_{0, \rm div} (\Omega )^n\Big)^{\circ }$ 
\,\footnotemark\,, is a closed subspace of $W^{ -1,\, q}(\Omega )^n$. 
\footnotetext{\,Let $X$ be a Banach space.  For a subset $M \subset X$ we define the annihilator $M^\circ$ 
which contains of all functionals $x^\ast\in X^\ast$ such that $\langle x^\ast, x \rangle=0$ 
for all $x\in  M$. Note, that $M^\circ$ is always a closed subspace of $X^\ast$.}
\end{rem}

\hspace*{0.5cm} 
Next, corresponding to a given functional $\bu^\ast\in G^{-1,\,q}(\Omega )^n $  
there exists a unique $[p] \in  L^q_{\rm pot}(\Omega )/\R $  such that 
\be\label{3.3}
\nabla_q p= \bu^\ast \quad \forall\, p\in  [p].
\ee
The element $[p]$  will be denoted by  $\overline{\mathscr{P}}_{q, \Omega }(\bu^ \ast)$, 
which defines a bijective linear mapping from $G^{-1,\,q}(\Omega )^n \rightarrow 
L^q_{\rm pot}(\Omega)/\R$.  

\begin{rem}
\label{rem3.5}
Let $1<q_1,q_2< +\infty $. From the above definition it is immediately clear that 
\be\label{3.4}
\overline{\mathscr{P}}_{q_1, \Omega } (\bu^\ast) =
\overline{\mathscr{P}}_{q_2, \Omega }(\bu^\ast )
\quad \forall\, \bu^\ast\in G^{-1,q_1}\cap G^{-1, q_2} (\Omega )^n.  
\ee
If no confusion can arise, we may omit  both subscript $q$ and  $\Omega $  and write 
$\overline{\mathscr{P}}$ in place of $\overline{\mathscr{P}}_{q, \Omega }$.  
\end{rem}

\hspace*{0.5cm} 
Let $\bu^\ast\in W^{ -1,\, q}(\Omega )^n$.  By $\nabla \cdot \bu^\ast$  we denote the functional
\[
v  \mapsto - \langle \bu^\ast, \nabla v \rangle,\quad v\in  W^{ 2,\, q'}_0(\Omega ),
\]
which belongs to $W^{ -2,\, q}(\Omega )$.  Then we have the following 

\begin{lem}
\label{lem3.6}
Let $\bu^\ast \in  G^{-1,\,q}(\Omega )^n $. Suppose, $\nabla \cdot \bu^\ast =0$. 
Then, every $p\in  \overline{\mathscr{P}}(\bu^\ast)$ is harmonic. 
\end{lem}

{\bf Proof}:  Note,  that $p\in  \overline{\mathscr{P}}(\bu^\ast)$ is equivalent to  
$p\in  L^q_{\rm pot}(\Omega ) $ with $\bu^\ast= \nabla _q p $.  Thus, in view of \eqref{3.1} we have 
\[
\intl_{\Omega } p \Delta \phi  dx = - \langle \nabla _q p, \nabla \phi  \rangle  
= - \langle \bu^\ast, \nabla \phi  \rangle  = \langle \nabla \cdot \bu^\ast, \phi  \rangle =0
\]
for all $\phi \in C^\infty _{\rm c}(\Omega ) $. Hence, by Weyl's lemma $p$ is harmonic. 
\hfill \Beweisende

\vspace{0.3cm}
\hspace*{0.5cm} 
Next,  we derive  some interesting properties of  functions $p \in L^q_{\rm pot} (\Omega )$.  
We begin with the following definition.

\begin{defin}
\label{def3.7}
A subdomain $U\subset \Omega $ is called $q$-{\it suitable} if 
$p|_{U} \in L^q(U)$ for all $p\in  L^q_{\rm pot}(\Omega ) $.  
\end{defin}

\begin{rem}
\label{rem3.8}
1. If $U \subset \Omega $ is $q$-suitable then $\mes U <+\infty $, since 
$1 \in L^q_{\rm pot}(\Omega )$.    

2. Every $U\Subset \Omega $ 
is $q$-suitable. 

3. If $ L^q_{\rm pot}(\Omega )=  L^q (\Omega )$ then every subdomain $U\subset \Omega $ 
is $q$-suitable.

\end{rem}

\vspace{0.2cm}
\hspace*{0.5cm} 
We have the following weak Poincar\'{e}-type inequality.

\begin{lem}[{\bf Weak Poincar\'{e}-type inequality}]
\label{lem3.9}
Let $U\subset \Omega $ be a $q$-suitable  subdomain. Then there exists a constant $c>0$ such that 
\be\label{3.5}
\|p-p_{U}\|_{L^q(U)}  \le  c  \|\nabla_q  p\|_{W^{ -1,\, q}(\Omega )}  
\quad \forall\, p \in L^{q}_{\rm pot}(\Omega ).  
\ee
\end{lem}

{\bf Proof}:  In view of Remark\,\ref{rem3.8} there holds $\mes U <+\infty $.  
On $G^{ -1,\, q}(\Omega )^n$ we introduce the following  equivalence relation. 
Let $\bu^\ast, \bv^\ast \in   G^{ -1,\, q}(\Omega )^n $. We say $\bu^\ast \sim_{U}  \bv^\ast$   if
\[
\bu^\ast|_{W^{ 1,\, q'}_0(U)} = \bv^\ast|_{W^{ 1,\, q'}_0(U)}. 
\]
We define the linear mapping $\Phi: G^{ -1,\, q}(\Omega )^n/\sim_{U}  \rightarrow L^q_0 (U)$ by setting $\Phi ([\bu^\ast]) : = p|_U - p_U $, where 
$p \in  \overline{\mathscr{P}}(\bu^\ast)$.  
By the assumption of the lemma, $\Phi $ is surjective. On the other hand, if $\Phi ([\bu^\ast])=0$  
there exists $p\in L^q_{\rm pot}(\Omega ) $  vanishing on $U$ such that 
$\bu^\ast =\nabla _q p $. Hence, $\langle \bu^\ast, \bu \rangle=0$ for all 
$\bu \in W^{ 1,\, q'}_0(U)^n$, which shows that $\Phi $ also is injective and hence bijective.  Its inverse is bounded, which follows from 
\begin{align*}
\|\Phi ^{-1} p\|_{W^{ -1,\, q}(\Omega )/\sim_{U} } &= \inf  \Big\{ \|\bu^\ast\|_{W^{ -1,\, q}}  \,\Big|\, \bu^\ast\in  \Phi ^{-1} p \Big\}
\\     
&\le \sup_{\substack{\bu \in W^{ 1,\, q'}_0(\Omega )
\\\|\bu\|_{W^{ 1,\, q'}\le1} } } 
\intl_{U} p \nabla \cdot \bu dx  \le  \|p\|_{L^q(U)} \|\nabla \cdot \bu\|_{L^q} 
\\
&\le n  \|p\|_{L^q(U)} 
\end{align*}
for all $p\in  L^q_0(U)$.  
By the closed range theorem we deduce that $\Phi $ also is bounded.    This implies for 
$p\in  L^q_{\rm pot}(\Omega )$ and $ \bu ^{ \ast}=\nabla_q( p-p_{ U})$
\begin{align*}
\|p-  p_{U} \|_{L^q(U)} &= \| \Phi ([ \bu ^{ \ast}])\|_{L^q(U)}
\\
& \le  c\|[ \bu ^{ \ast}]\|_{W^{ -1,\, q}(\Omega )/ \sim_{U}  }        
\\     
&\le c \|\nabla_q (p-p|_{U}) \|_{W^{ -1,\, q}} =  c \|\nabla_q p \|_{W^{ -1,\, q}}. 
\end{align*}
This proves  \eqref{3.5}. \hfill \Beweisende

\vspace{0.3cm}
\hspace*{0.5cm} 
In what follows, let   $U\subset \Omega $ be a fixed $q$-suitable  domain. 
Then, for $\bu^\ast\in  
G^{-1,\,q}(\Omega )^n$ by $ \mathscr{P}^{(U)} (\bu^\ast)$  we denote the unique pressure 
$p \in  \overline{\mathscr{P}} (\bu^\ast )$ fulfilling $p_U=0$.  
If $L^q_{\rm pot} (\Omega)= L^q(\Omega )$ we might take $U=\Omega$.  
In this case, we shortly write $\mathscr{P} $ in place of $\mathscr{P}^{(\Omega )} $.

\vspace{0.2cm}
\hspace*{0.5cm} 
From Lemma\,\ref{lem3.9} we easily derive the following 

\begin{lem}
\label{lem3.10}
For every $q$-suitable subdomain $G\subset \Omega $  there exists a constant $c_{G}>0 $ such that 
\be\label{3.6}
\|\mathscr{P}^{(U)} (\bu^\ast) \|_{L^q(G)} \le c_G \|\bu^\ast\|_{W^{ -1,\, q}(\Omega )}  
\quad \forall\, \bu^\ast\in  G^{-1,\,q}(\Omega )^n.  
\ee
\end{lem}

{\bf Proof}:  Firstly, assume that $G\cap  U\neq \emptyset $.  Let  $\bu^\ast \in 
G^{- 1,\, q}(\Omega )^n$,  and set $p :=  \mathscr{P}^{(U)} (\bu^\ast) $. 
Clearly, as $p_{U}=0 $ we easily  find 
\begin{align}
&\frac {1} {\mes G } \intl_{G} |p|^q  dx 
\cr
&\qquad\le     \frac {1} {\mes G } \intl_{G} |p- p_G|^q  dx 
+ 4^{s-1}  |p_{G}- p_{G\cap U}|^q+  4^{s-1} | p_{G\cap U} - p_U|^q.  
\label{3.7}
\end{align} 
Estimating the second  term on the right of \eqref{3.7}  by means of 
\[
|p_{G}- p_{G\cap U}|^q \le \frac {2^q} {\mes (G\cap U)} \intl_{G} |p- p_{G} |^q dx,
\]
and the third term by a similar one, we are led to 
\[
\intl_{G} |p|^q  dx \le c \intl_{G} |p- p_{G} |^q dx + c\intl_{U} |p- p_{U} |^q dx. 
\]
Recalling that both $U$  and $G$ are $q$-suitable, we are in a position to apply Lemma\,\ref{lem3.9} to both terms on the   right-hand  side of the above estimate. This implies  \eqref{3.6}. 

\hspace*{0.5cm} 
Secondly, if $G \cap  U=\emptyset $, we take $G\subset G_0 \subset \Omega$ such that 
$G_0$ is $q$-suitable and $ G_0\cap U \not= \emptyset$.  Then \eqref{3.6} immediately follows from the first case.
\hfill \Beweisende 

\begin{rem}
\label{rem3.11}
1. Lemma\,\ref{lem3.10} says that for every $q$-suitable $G\subset \Omega $, the mapping 
$\bu^\ast \mapsto \mathscr{P}^{(U)}(\bu^\ast)|_{G} $ 
is a bounded linear operator from $G^{-1,\,q} (\Omega )^n$  into $L^q(G)$. 

\vspace{0.2cm}
2. If $L^q_{\rm pot} (\Omega )=L^q(\Omega ) $ then $\mathscr{P}: G^{-1,\,q}(\Omega )^n \rightarrow L^q(\Omega )$ is bounded. 

\end{rem}

%%%%%%%%%%%% 
% 
%section 4 
% 
%%%%%%%%%%%% 
\section{Projections onto $G^{-1,\,q}(\Omega )^n$} 
\setcounter{secnum}{\value {section}
\setcounter{equation}{0}
\renewcommand{\theequation}{\mbox{\arabic{secnum}.\arabic{equation}}}}

\hspace*{0.5cm} 
As we have  seen in the previous section functionls in $W^{ -1,\, q}(\Omega )^n$   vanishing 
on $C^\infty _{\rm c, div}(\Omega )^n $ equal to functionals of the form $\nabla_q p $ with potential $p$.  
Furthermore,  in many applications such functionals are expressed by a sum, namely  
$\nabla _q p =\bu_1^\ast + \bu_2^\ast$ such that $\bu_i^\ast\in W^{- 1,\, q_i}(\Omega )^n$ \, $ (1<q_i<+\infty)$
\, $(i=1,2)$.  If  there is  an operator $\bE^\ast$
which simultaneously projects $\bu_1^\ast$ into $G^{-1, q_1}(\Omega )^n $ and 
$\bu_2^\ast$ into $G^{-1, q_2}(\Omega )^n $ we are able to write 
\be\label{4.1}
\nabla_q p = \bE^\ast(\nabla_q p)= 
\nabla _{q_1} p_1 + \nabla _{q_2} p_{2},\quad p_1\in L^{q_1}_{\rm pot}(\Omega ), 
\quad p_2\in L^{q_2}_{\rm pot}(\Omega).        
\ee
Unfortunately,  the existence of an operator $\bE^\ast$ which implies  \eqref{4.1} 
isn't necessarily guaranteed,  unless  the domain $\Omega $ enjoys certain regularity properties. 
On the other hand, if such projection $\bE^\ast$  exists,  there are infinite projections 
leading to \eqref{4.1}. 
Nevertheless,  among all such projections  there is a canonical one which is related to the 
existence and uniqueness of weak solutions to the Stokes system (or a Stokes-like system for unbounded domains). 

\hspace*{0.5cm}       
Unlike the case $q\neq 2$, the case $q=2$ appears to be special, due to the Hilbert space structure we are in a position  to define $\bE^\ast_2$ as the orthogonal projection of $W^{ -1,\, 2}(\Omega )^n$ 
onto the closed subspace $G^{-1, \,2}(\Omega )^n $, where the scalar product in $W^{ -1,\, 2}(\Omega )^n$ 
is given  by 
\be\label{4.2}
((\bu^\ast, \bv^\ast))_\ast = \langle \bu^\ast, \bJ_2^{-1} \bv^\ast  \rangle = 
((\bJ_2^{-1} \bu^\ast, \bJ_2^{-1} \bv^\ast)),\quad \bu^\ast, \bv^\ast \in  W^{ -1,\, 2}(\Omega )^n.
\ee
Here  $((\cdot, \cdot ))$ denotes the usual scalar product in $W^{ 1,\, 2}_0(\Omega )^n$.  
Then, by $\bE_2$ we denote the dual of $\bE^\ast_2$ which appears to be a projection from $W^{ 1,\, 2}_0(\Omega )^n$ onto a closed subspace.  
Then then  \eqref{1.2} and \eqref{1.3} are fulfilled.  
In fact, \eqref{1.2} follows by the aid of the closed range theorem together with Lemma\,\ref{lem3.3}, as 
\[
{\rm im}\,\bE_2^\ast = G^{ -1,\, 2}(\Omega )^n = (W^{ 1,\, 2}_{0, \rm div} (\Omega ))^\circ 
= ({\rm ker}\, \bE_2)^\circ.  
\] 
In addition, there holds 
\be\label{4.3}
\bE^\ast_2 \bJ_2 = \bJ_2 \bE_2 \,\footnotemark\,. 
\ee
 \footnotetext{\, This can be readily seen by $\langle \bE^\ast_2 \bJ_2\bu, \bv \rangle = 
((\bE^\ast_2 \bJ_2\bu, \bJ_2\bv))_\ast = ((\bJ_2\bu, \bE^\ast_2\bJ_2\bv))_\ast=
 \langle \bJ_2 \bv, \bE_2 \bu \rangle =\langle \bJ_2 \bE_2 \bu, \bv \rangle $  for all 
$\bv\in  W^{ 1,\,2}_0(\Omega )^n$.}
Noting  that $\bJ_2 (\nabla \phi ) \in G^{-1,\,2}(\Omega )^n $  for all 
$\phi \in C^\infty _{\rm c}(\Omega ) $, the property \eqref{4.3} 
furnishes \eqref{1.3}, i.\,e.
\be\label{4.4}
\bE_2 (\nabla \phi ) = \nabla \phi \quad \forall\, \phi \in C^\infty _{\rm c}(\Omega ).  
\ee

\hspace*{0.5cm} 
Bearing in mind \eqref{4.3}, we give the   following definition of the  projection $\bE_{q'} $ for $1<q<+\infty $.  

\begin{defin}
\label{def4.1}
Let $1<q<+\infty $.  A domain $\Omega \subset \R^n$ is called  $\nabla _q $-{\it regular}  
if there exists a projection $\bE_{q'} : W^{ 1,\, q'}_{0} (\Omega )^n
\rightarrow W^{ 1,\, q'}_{0 } (\Omega )^n$ fulfilling 
\begin{align}
{\rm ker}\, \bE_{q'}  &= W^{ 1,\, q'}_{0, \rm div} (\Omega )^n,
\label{4.5}  
\\
\bE_{q'}(\bfpsi ) &= \bE_2 (\bfpsi) \quad \forall\, \bfpsi\in C^\infty _{\rm c}(\Omega )^n. 
\label{4.6}
\end{align}
\end{defin}

\begin{rem}
\label{rem4.2}
1. As $C^\infty_{\rm c} (\Omega )^n $ is dense in 
$W^{ 1,\, q'}_0(\Omega )^n$  the projection $\bE_{q'} $ 
is uniquely defined by \eqref{4.6}.  

\vspace{0.2cm}
2.  From \eqref{4.6}, by using \eqref{4.4} we immediately get 
\be\label{4.7}
E_{q'}(\nabla \phi ) = \nabla \phi \quad \forall\, \phi \in C^\infty _{\rm c}(\Omega ).   
\ee
 
3. If $W^{ 1,\, q'}_{0, \rm div} (\Omega )^n
= W^{ 1,\, 2}_{0, \rm div} (\Omega )^n \cap  W^{ 1,\, q'}_0(\Omega )^n$  then 
\eqref{4.5}  follows  from \eqref{4.6} since   \eqref{4.6} implies  
$\bfpsi - \bE_{q'}\bfpsi =\bfpsi - \bE_{2}\bfpsi \in  
W^{ 1,\, 2}_{0, \rm div} (\Omega )^n \cap  W^{ 1,\, q'}_0(\Omega )^n $  for all $\bfpsi \in  C^\infty _{\rm c}(\Omega )^n $.  

\vspace{0.2cm}
4.  The dual operator  $\bE^\ast_q=(\bE_{ q'})' : W^{ -1,\, q}(\Omega ) ^n\rightarrow 
W^{ -1,\, q}(\Omega )^n$  defines  a projection onto $G^{-1,\,q}(\Omega )^n $.  
By the aid of the closed range theorem we infer 
\begin{align}
W^{ -1,\, q}(\Omega )^n = G^{-1,\,q}(\Omega )^n \oplus  {\rm ker} \,\bE^\ast_q, 
\label{4.8}
\\
W^{1,\, q'}_0(\Omega )^n = {\rm im} \,\bE_{q'}\oplus W^{ 1,\, q}_{0, \rm div} (\Omega )^n.
\label{4.9}
\end{align}
If $\bJ_q$ is an isomorphism then 
\begin{align}
{\rm im} \,\bE_{q'} &= \{ \bv\in  W^{ 1,\, q'}_0(\Omega )^n \,|\, \bJ_{ q'}\bv \in G^{ -1,\, q'} (\Omega )^n\},
\label{4.10}
\\
{\rm ker} \,\bE^\ast_q &= \{ \bJ_q \bv \,|\,\bv \in W^{ 1,\, q}_{0, \rm div} (\Omega )^n\}.
\label{4.11}
\end{align}
\end{rem}

\begin{rem}
\label{rem4.3}
Let $1<q_1,q_2 <+\infty $.  Suppose $\Omega $ is  $\nabla _{q_1} $-regular and also  
$\nabla _{q_2} $-regular. 
As $C^\infty_{\rm c }(\Omega ) ^n $  
is dense in both $W^{ 1,\, q_1}_0(\Omega )^n$ and $W^{ 1,\, q_2}_0(\Omega )^n$, 
\eqref{4.6} implies $\bE_{q_1} ^\ast \bu^\ast = \bE_{q_2} ^\ast \bu^\ast $  for all 
$\bu^\ast \in W^{ -1,\, q_1}\cap  W^{ 1,\, q_2}(\Omega )^n$.   
Thus, we may write shortly $\bE^\ast$  in place of $\bE_{q_1} ^\ast$ or $\bE_{q_2} ^\ast$. 
In particular, $\bE^\ast $ is a projection in $ W^{ -1,\, q_1}\cap  W^{ -1,\, q_2}(\Omega )^n$.    

\hspace*{0.5cm} 
On the other hand,  let $ \bu^\ast =  \bu^\ast _1+\bu^\ast_2 \in  
(W^{ -1,\, q_1}+ W^{ -1,\, q_2} )(\Omega )^n$, then 
$\bE^\ast\bu^\ast  = \bE^\ast_{q_1}  \bu_1 ^\ast + \bE_{q_2}  \bu^\ast_2
\in  (W^{ -1,\, q_1}+ W^{ -1,\, q_2} )(\Omega )^n$.  Elementary, 
\[
\begin{cases}
\langle \bE^\ast\bu^\ast, \bv \rangle \le c\max \{\|\bu_1^\ast\|_{W^{- 1,\, q_1}} , 
\|\bu_2^\ast \|_{W^{ -1,\, q_2}} \} (\|\bv\|_{W^{ 1,\, q_1}}+  \|\bv\|_{W^{ 1,\, q_2}})
\\[0.2cm]
\forall\,  \bv\in  W^{ 1,\, q_1}\cap  W^{ 1,\, q_2}(\Omega )^n.
\end{cases}
\]
Hence, $\bE^\ast$ is a projection in $ (W^{ -1,\, q_1}+  W^{ -1,\, q_2})(\Omega )^n $ which fulfills  
\eqref{4.1}.  The dual of $\widetilde{\bE}$ of $\bE^\ast$  is a projection from  
$ W^{ 1,\, q_1'}_0\cap  W^{ 1,\, q_2'}_0(\Omega )^n $ into itself satisfying 
$\widetilde{\bE} \bfpsi = \bE_2 \bfpsi $  for all 
$\bfpsi \in C^\infty _{\rm c } (\Omega )^n$. However, we don't know whether 
$\widetilde{\bE} \bv = \bE_{q_1}  \bv= \bE_{q_2}  \bv $ for all 
$ \bv \in W^{ 1,\, q_1}_{0} \cap W^{ 1,\, q_2}_{0} (\Omega )^n$.  Eventually, this property 
holds if    $ C^\infty _{\rm c } (\Omega )^n $ is dense in 
$W^{ 1,\, q_1}_{0} \cap W^{ 1,\, q_2}_{0} (\Omega )^n$ which is true for uniform Lipschitz domains. 

\end{rem}

\vspace{0.1cm}
\hspace*{0.5cm} 
The condition \eqref{4.6}  implies the following properties of $\bE_{q'}. $

\begin{lem}
\label{lem4.4}
Let $1<q<+\infty $, and let $\Omega $  be $\nabla_q$-regular. Then the following statements are true.

\vspace{0.2cm}
1.  There holds 
\be\label{4.12}
\bJ_{q'} \bE_{q'}(\bfpsi ) = \bE_q^\ast \bJ_q (\bfpsi )\quad \forall\, \bfpsi\in  C^\infty _{\rm c}(\Omega )^n.   
\ee

2. Let $\bu^\ast\in  W^{ -1,\, q}(\Omega )^n$ such that $\langle \bu^\ast, \nabla \phi  \rangle=0$ 
for  all $\phi \in C^\infty _{\rm c}(\Omega ) $. Then,  every potential 
$p \in L^q_{\rm pot }(\Omega )$ of $\bE^\ast_q \bu^\ast$  is harmonic. 
\end{lem}

{\bf Proof}: 1.  According to \eqref{4.6} and \eqref{4.3},  we get 
for every $\bfpsi, \bfphi  \in  C^\infty _{\rm c}(\Omega )^n $
\begin{align*}
\langle \bE_q^\ast \bJ_q (\bfpsi ), \bfphi  \rangle & =\langle \bJ_{q} \bfpsi, \bE_{q'} (\bfphi )    \rangle 
= \langle \bfpsi, \bJ_{2} \bE_2 (\bfphi) \rangle 
\\     
&= \langle \bfpsi, \bE ^\ast_2\bJ_{2} (\bfphi) \rangle
\langle \bJ_2\bE_2(\bfpsi), \bfphi \rangle
= \langle \bJ_{q'} \bE_{q'} (\bfpsi), \bfphi \rangle.
\end{align*}
Whence, \eqref{4.12}.

2. Let $p\in  L^q_{\rm pot}(\Omega ) $ such that $\nabla _q p = \bE^\ast_q \bu^\ast$. 
In light of \eqref{4.7} we calculate 
\[
\intl_{\Omega } p \Delta \phi  dx = - \langle \bE^\ast_q \bu^\ast, \nabla \phi  \rangle  
= - \langle \bu^\ast, \nabla \phi  \rangle  = 0
\]
for all $\phi \in C^\infty _{\rm c}(\Omega ) $. Hence, by Weyl's lemma $p$ is harmonic. 
\hfill \Beweisende

\vspace{0.2cm}
\hspace*{0.5cm} 
Next, we turn to the decomposition of the pressure by using the projection $\bE^\ast$.

\begin{thm}
\label{thm4.5}
Let $1<q_i< +\infty $ \, $(i=1,\ldots,N)$. Suppose $\Omega $ is  $\nabla _{q_i}  $-regular 
for all $i=1,\ldots,N$.  Let $\bu^\ast_i \in  
W^{ -1,\, q_i}(\Omega )^n $  \, $(i=1,\ldots,N)$ such that 
\be\label{4.13}
\suml_{i=1}^{N} \langle \bu^\ast_i, \bfpsi  \rangle = 0\quad 
\forall\, \bfpsi \in  \bigcap_{i=1}^N W^{ 1,\, q_i'}_{0, \rm div}  (\Omega )^n. 
\ee
Then, for every $U\Subset \Omega $  there holds
\be\label{4.14}
\suml_{i=1}^{N} \langle \bu^\ast_i, \bfpsi  \rangle = -\suml_{i=1}^{N} \intl_{\Omega} p_i \nabla  \cdot \bfpsi dx  
\quad \forall\, \bfpsi\in C^\infty _{\rm c}(\Omega )^n,  
\ee
where $p_i = \mathscr{P}^{(U)} (\bE^\ast_{q_i} \bu^\ast_{i})$ \, $ (i=1,\ldots,N)$.      

\end{thm}

{\bf Proof}:   Let $\bfpsi \in  C^\infty _{\rm c}(\Omega )^n $ be arbitrarily chosen.  
Observing \eqref{4.6}, we get $\bfpsi - \bE_{2}(\bfpsi ) = \bfpsi - \bE_{q_i}(\bfpsi) \in  
W^{ 1,\, q_i}_{0, \rm div} (\Omega )^n$  for all $i=1,\ldots,N$, in view of \eqref{4.13} we infer  
\[
\suml_{i=1}^{N} \langle \bE_{q_i}^\ast \bu_i^\ast, \bfpsi   \rangle
= \suml_{i=1}^{N} \langle \bu_i^\ast, \bE_{q_i'} \bfpsi   \rangle 
= \suml_{i=1}^{N} \langle \bu_i^\ast, \bfpsi   \rangle. 
\]
Let $U\Subset \Omega $. According to Remark\,\ref{rem3.8}, $U$ is $q$-suitable. Setting 
$p_i =\mathscr{P}^{(U)} (\bE^\ast_{q_i} \bu^\ast_{i})$ \, $ (i=1,\ldots,N)$, as $ \bu ^{ \ast}_i= \nabla _{ q_i} p_i$,  the identity 
\eqref{4.14} follows from the latter identity.   This completes the proof of the theorem. 
\hfill \Beweisende   

\begin{rem}
\label{rem4.6}
1. Let  $\Omega \subset \R^n $ be a $\nabla _q$-regular domain $(1<q<+\infty )$, and 
let $U\Subset \Omega $ be $q$-suitable. Then for every $\bu^\ast\in 
W^{ -1,\, q}(\Omega )^n$ there exists a unique associate pressure $p \in L^q_{\rm pot}(\Omega ) $
with $p_{U}=0 $, defined by $p:= \mathscr{P}^{(U)} (\bE_q^\ast\bu^\ast) $.  This pressure will be 
denoted shortly by $\mathscr{P}^{(U)} (\bu^\ast)$. In fact,  $\mathscr{P}^{(U)}$ defines 
a linear map from $ W^{ -1,\, q}(\Omega )^n $ onto $L^q_{\rm pot}(\Omega ) $. In particular, there holds 
\be\label{4.15}
\langle \bE_q \bu^\ast, \bfpsi  \rangle = 
-\intl_{\Omega} \mathscr{P}^{(U)} (\bu^\ast) 
\nabla \cdot \bfpsi  dx  \quad \forall\, \bfpsi \in C^\infty_{\rm c}(\Omega ) ^n.   
\ee
On the other hand for every subdomain $G\subset \Omega $ being $q$-suitable, the mapping 
$\bu^\ast \mapsto \mathscr{P}^{(U)} (\bu^\ast) |_{G} $ is a bounded linear operator from 
$ W^{ -1,\, q}(\Omega )^n $ onto $L^q (G)$. 

\hspace{0.5cm}
Respectively, if 
$L^q_{\rm pot} (\Omega )= L^q(\Omega )$, we define $\mathscr{P} (\bu^\ast) := \mathscr{P}(\bE^\ast_q \bu^\ast)$ such that  $\mathscr{P} $ is a bounded linear operator from 
$ W^{ -1,\, q}(\Omega )^n $ onto $L^q (\Omega )$ (see Remark\,\ref{rem3.11} for the definition of $ \mathscr{P}(\bE^\ast_q \bu^\ast)$ ). 
\end{rem}

%%%%%%%%%%%% 
% 
%section 5 
% 
%%%%%%%%%%%% 
\section{Sufficient condition for $\nabla _q$-regularity} 
\setcounter{secnum}{\value {section}
\setcounter{equation}{0}
\renewcommand{\theequation}{\mbox{\arabic{secnum}.\arabic{equation}}}}

In this section we wish to present a sufficient condition for $\Omega $ being $\nabla _q$-regular, 
based on the existence and uniqueness of $q$-weak solutions to the system 
\begin{align}
  -\Delta \bu + \delta \bu &= \bu^\ast -\nabla p,\quad \nabla \cdot \bu =0\quad \mbox{in}\quad 
\Omega,
\label{5.1} 
\\
\bu &= {\bf 0} \quad \mbox{on}\quad \partial \Omega,
\label{5.2}
\end{align}
with $\delta =0$ or $\delta =1$. If $\delta =0$,  \eqref{5.1}, \eqref{5.2}  forms the Stokes system. In case  $\delta =1$ we call \eqref{5.1}, \eqref{5.2}  the  {\it Stokes-like system}. Concerning weak solutions to  the Stokes or the Stokes-like system we give 
the following definition.

\begin{defin}
\label{def5.1}
Let $1<q<+\infty $. Let $\bu^\ast \in  W^{ -1,\, q}(\Omega )^n$. 
Then $\bu\in W^{ 1,\, q}_{0, \rm div} (\Omega )^n$  is called a {\it $q$-weak solution} to \eqref{5.1}, \eqref{5.2} 
if $\nabla \cdot \bu=0$ a.\,e. in $\Omega $, and there holds 
\be\label{5.3}
\intl_{\Omega} \nabla \bu: \nabla \bv + \delta \bu \cdot \bv dx = \langle \bu^\ast, \bv \rangle 
\quad \forall\, \bv \in  W^{ 1,\, q'}_{0, \rm div} (\Omega )^n.
\ee
\end{defin}

\begin{rem} 
\label{rem5.2}
1. Using the canonical embedding $L^q(\Omega )^n \hookrightarrow 
W^{ -1,\, q}(\Omega )^n$,  and employing  Lemma\,\ref{lem3.3}, 
we see that $\bu\in  W^{ 1,\, q}_{0, \rm div} (\Omega )^n $ is a $q$-weak 
solution to \eqref{5.1}, \eqref{5.2} iff there exists $p\in  L^q_{\rm pot}(\Omega ) $ such that  
\be\label{5.4}
\bu^\ast -\Delta_q  \bu -\delta \bu= -\nabla_q p \quad  \mbox{in} \,\, W^{-1,\,q}(\Omega )^n,
\ee
where 
$  \langle \Delta_q \bu, \bv  \rangle = \intl_{\Omega } \nabla \bu : \nabla \bv   dx$, 
$ \bu \in W^{1,\, q}_{0, \rm div}(\Omega )^n, \bv \in W^{1,\, q'}_{0, \rm div}(\Omega )^n$.

\vspace{0.2cm}
2.  Clearly, \eqref{5.3}  can be interpreted as an operator equation
\be\label{5.5}
\bT_q \bu = \bu^\ast|_{W^{ 1,\, q'}_{0, \rm div} (\Omega )}, 
\ee
where $\bT_q \bu$ stands for the restriction of $-\Delta _q \bu + \delta \bu $ to 
$ W^{ 1,\, q'}_{0, \rm div} (\Omega )^n$ which appears to be a bounded  linear operator 
from $W^{ 1,\, q}_{0, \rm div} (\Omega )^n$ into $(W^{ 1,\, q'}_{0, \rm div} (\Omega )^n)'$.  
By using a routine  functional analytic argument, we infer that  the existence and uniqueness of 
$q$-weak solutions to the system \eqref{5.1}, \eqref{5.2}  
is equivalent to   $\bT_q$ being an isomorphism.  Furthermore, having  
\be\label{5.6}
(\bT_q)' = \bT_{q'}\quad \forall\, 1<q<+\infty,  
\ee
the existence and uniqueness of $q$-weak solutions to \eqref{5.1}, \eqref{5.2} is equivalent to the existence and uniqueness  of $q'$-weak solutions  
to \eqref{5.1}, \eqref{5.2}.  

\end{rem}
 
\hspace*{0.5cm} 
Next we shall introduce a sufficient condition  for $\Omega $ being $\nabla _q$-regular.  

\begin{lem}
\label{lem5.3}
Let $\Omega \subset \R^n$ be a domain. Let $1<q<+\infty$.  Suppose that  

\begin{itemize}
\item[(a)] For every 
$\bu^\ast\in W^{ -1,\, q}(\Omega )^n$ there exists a unique $q$-weak solution 
$\bu\in W^{ 1,\, q}_{0, \rm div} (\Omega )^n$ to \eqref{5.1}, \eqref{5.2}
$(\delta =0$  if $\Omega $ is bounded or    $\delta =1$ otherwise$)$.  

\item[(b)]
For every $\bfpsi \in C^\infty _{\rm c} (\Omega )^n$  the $q$-weak solution to \eqref{5.1}, \eqref{5.2} 
with  $ \bu ^{ \ast}=- \Delta \bfpsi  +\delta \bfpsi $ belongs to  $W^{ 1,\, 2}_{0, \rm div} (\Omega )^n$.  
\end{itemize}

Then $\Omega $ is $\nabla_{q} $-regular. 
\end{lem}

{\bf Proof}:  1.  Let $\bv \in W^{ 1,\, q'}_0(\Omega )^n$.  Note that, according to 
Remark\,\ref{rem5.2}/2.,  (a) continuous to hold after replacing $q$ by $q'$ therein. 
Thus, there exists a unique $q'$-weak solution $\bu\in  W^{ 1,\, q'}_{0, \rm div} (\Omega )^n$  to 
\eqref{5.1}, \eqref{5.2} with  right-hand  side $\bu^\ast = - \Delta _{q'}  \bv + \delta \bv
= \bJ_{q'}  \bv$.  Setting $\bE_{q'}  \bv =\bv - \bu$, and   having  $\|\bE_{q'}  \bv \|_{W^{ 1,\, q'}} \le 
c \|\bv\|_{W^{ 1,\, q'}} $,  we see that $\bE_{q'} $ is a linear bounded operator from 
$W^{ 1,\, q'}_0(\Omega )^n$ into itself.  
By the definition of $\bE_{q'} $ it is readily seen that ${\rm ker}\, \bE_{q'}  = W^{ 1,\, q'}_{0, \rm div} (\Omega )^n$.  
In particular, owing to  $\bv- \bE_{q'}  \bv \in  W^{ 1,\, q'}_{0, \rm div} (\Omega )^n$  it follows 
$\bE_{q'} (\bv - \bE_{q'}  \bv)={\bf 0}$ for all $\bv \in W^{ 1,\, q'}_{0 } (\Omega )^n $. 
Consequently, $\bE^2_{q'}  = \bE_q$. This, shows that $\bE_{q'} $ is a projection enjoying \eqref{4.5}. 

\hspace*{0.5cm} 
2. Now, it remains to verify the condition \eqref{4.6}.  To see this, let $\bfpsi \in C^\infty _{\rm c}(\Omega )^n $ 
be arbitrarily chosen. Let $\bu\in  W^{ 1,\, q'}_{0, \rm div} (\Omega )$ denote the unique 
$q'$-weak solution to \eqref{5.1}, \eqref{5.2}  with right-hand side $\bu^\ast =-\Delta  \bfpsi +\delta \bfpsi $.  
Due to (b) we have  $\bu \in  W^{ 1,\, 2}_{0, \rm div} (\Omega )^n$. Recalling the definition 
of $\bE_2$, we get $\bE_2 \bfpsi = \bfpsi - \bu = \bE_{q'}  \bfpsi $. Whence, \eqref{4.6}. This shows that 
$\Omega $ is $\nabla _{q} $-regular. 
\hfill \Beweisende

\begin{rem}
\label{rem5.4}
1.  Recalling the definition of  $\bE_2 $,  we see that $I-\bE_2$ becomes the orthogonal projection 
onto $W^{ 1,\, 2}_{0, \rm div} (\Omega )^n$.  Let $\bv  \in  W^{ 1,\, 2}_0(\Omega )^n$. Define, 
$\bu = \bv - \bE_2 \bv$  in view of \eqref{4.3}  we see that 
\[
\bJ_2 \bu = -\Delta _2 \bu + \delta \bu = \bJ_2 \bv - \bE_2^\ast \bJ_2 \bv. 
\]
Thus, $\bv - \bE_2 \bv$ is the unique $2$-weak solution to \eqref{5.1}, \eqref{5.2} with right-hand side $\bu^\ast= \bJ_2 \bv$.  On the other hand for any  $\bu^\ast $, the function  
$\bu= \bJ^{-1}_2\bu^\ast - \bE_2^\ast\bu^\ast  $  becomes the unique $2$-weak solution to 
\eqref{5.1}, \eqref{5.2}. 

\vspace{0.2cm}
2. If $\Omega \subset \R^n$  is bounded the condition (a) implies that $\Omega $ is $\nabla _q$-regular and $\nabla _{q'} $-regular as well.  

\vspace{0.2cm}
\hspace*{0.5cm} 
Indeed, for $1<q <2$  verifying  
$W^{ 1,\, q'}_{0, \rm div} (\Omega )^n \subset 
W^{ 1,\, 2}_{0, \rm div} (\Omega )^n$ we immediately get $\bu =\bfpsi - \bE_{q'} \bfpsi  \in W^{1,\, 2}_{0, \rm div} (\Omega )^n$,  
and thus (b) holds.  If $2\le q<+\infty $ we argue as follows.   
Let $\bfpsi\in C^\infty _{\rm c} (\Omega )^n $. As we have mentioned above in  1., 
$\bu = \bfpsi - \bE_2\bfpsi \in W^{ 1,\, 2}_{0, \rm div} (\Omega )^n$ is the unique $2$-weak solution to \eqref{5.1}, \eqref{5.2} with $\bu^\ast= -\Delta \bfpsi $. 
Since  $W^{ 1,\, 2}_{0, \rm div} (\Omega )^n \subset 
W^{ 1,\, q'}_{0, \rm div} (\Omega )^n$,  from (a)  it follows that 
$\bu$ is the $q'$-weak solution to the corresponding system.  This shows that $(b)$ holds. According to Lemma\,\ref{lem5.3} $\Omega $ is $\nabla _{q'} $-regular, and in view of Remark\,\ref{rem5.2}, $ \Omega $ is  $\nabla _{q}$-regular.

\end{rem}

\vspace{0.2cm}
\hspace*{0.5cm} 
Next, we shall  present examples of $\nabla _q$-domains, which 
occur in various applications. However, the list of domains below is not  complete, and 
it will be left to the reader to find more relevant examples used in further applications.

\subsection{The whole space $\R^n$} 

By employing the well-known Calder\'{o}n-Zygmund inequality along with the classical regularity theory 
and Sobolev's inequalities, one verifies that the operator 
$-\Delta + I: W^{ 2,\, q}(\R^n) \rightarrow L^q(\R^n)$ is an isomorphism for all $1<q< +\infty$.  
By using an interpolation argument, the above statement implies that $\bJ_q= - \Delta _q + \bI: 
W^{ 1,\, q}(\R^n)^n \rightarrow W^{ -1,\, q}(\R^n)^n$ becomes an isomorphism for all $1<q<+\infty $.    
This yields the following 

\begin{thm}
\label{thm5.5}
Let $1<q<+\infty $.  For every $\bu^\ast \in W^{ -1,\, q}(\R^n)^n$ there exists a unique 
$q$-weak solution $\bu\in W^{ 1,\, q}_{0, \rm div} (\R^n)^n$ to \eqref{5.1}, \eqref{5.2} 
for  $\delta =1$. 
\end{thm}

{\bf Proof}: First, note that $W^{ 1,\, q}_{0, \rm div} (\R^n)^n$ equals the space of all 
$\bu \in W^{ 1,\, q}(\R^n)^n$ with $\nabla \cdot \bu=0$.  This easily follows from 
\[
W^{ 1,\, q}(\R^n)^n = W^{ 1,\, q}_{0, \rm div} (\R^n)^n\oplus \{  \nabla f\,|\,f\in W^{ 2,\, q}(\R^n)\}
\]
which can be proved by the aid of Calder\'{o}n-Zymund's inequality  along with a duality argument.  
Now, let $\bP_q  $ denote the usual Helmholtz projection defined by 
\[
\bP_q\bu = \bu - \nabla  \Delta ^{-1} \nabla \cdot \bu, \quad \bu \in W^{ 1,\, q}(\R^n)^n.  
\]
Again using Calder\'{o}n-Zymund's inequality we see that $\bP_q $ is a projection  operator 
from $ W^{ 1,\, q}(\R^n)^n $ onto $ W^{ 1,\, q}_{0, \rm div} (\R^n)^n $.  

\hspace*{0.5cm} 
Now, let $\bu^\ast\in  W^{ -1,\, q}(\R^n)^n$ be arbitrarily chosen.  Set $\bu = \bJ^{-1} _{q} \bP_{q'}^\ast  \bu^\ast$. 
Verifying $\bJ^{-1} _{q} \bP_{q}^\ast = \bP_{q} \bJ^{-1} _{q'}$, we see that 
$\bu  \in  W^{ 1,\, q}_{0, \rm div} (\R^n)^n$ is a $q'$-weak solution to \eqref{5.1}, \eqref{5.2}.  
This solution also is unique, since $\bJ_{q} \bw \in G^{-1,\,q}(\Omega )^n$ implies 
$\bP^\ast_q\bJ_{q} \bw  = \bJ_{q} \bP_q \bw= \bJ_{q} \bw ={\bf 0}$, and thus $\bw={\bf 0}$. 
Whence, the assertion of the  theorem is proved. \hfill \Beweisende

\vspace{0.3cm}
\hspace*{0.5cm} 
As Theorem\,\ref{thm5.5} shows, condition (a) of Lemma\,\ref{lem5.3} is satisfied. Furthermore, as
$\bP_{q} \bfpsi = \bP_2 \bfpsi $  for all $\bfpsi \in C^\infty _{\rm c}(\Omega )^n $, condition 
(b)  also is fulfilled. Thus, by Lemma\,\ref{lem5.3}  we immediately get 

\begin{cor}
\label{cor5.6}
$\R^n$ is $\nabla _{q} $-regular for all $1<q< +\infty $. 
\end{cor}

\subsection{Bounded domains} 

To apply Lemma\,\ref{lem5.3} for bounded domains it will be sufficient to recall the existence and uniqueness of weak solutions to the Stokes system, which has been proved first by Cattabriga in  
\cite{Catt1961} for a bounded three dimensional $ C^2$ domain. 
For the  general case we  quote from \cite{GaSiSo1994} the following

\begin{thm}[{\bf Galdi, Simader, Sohr}] 
\label{thm5.7}
Let $\Omega \subset \R^n$ be a bounded $C^1$-domain. 
Let $1<q< +\infty$.  For every $\bu ^\ast\in  W^{ -1,\, q}(\Omega )^n$ and $g\in  L^q_0(\Omega )$  
there exits a unique pair $(\bu, p) \in  W^{ 1,\, q}_{0} (\Omega )^n\times L ^q_0(\Omega )$  
such that 
\be\label{5.7}
\begin{cases}
&\nabla \cdot \bu =g\quad \mbox{a.\,e. in}\quad \Omega 
\\
&- \Delta_q \bu = \bu^\ast - \nabla_q p \quad \mbox{in}\quad W^{ -1,\, q}(\Omega )^n. 
\end{cases}
\ee
In addition there holds
\be\label{5.8}
\|\nabla \bu\|_{\bL^q(\Omega )}+ \|p\|_{L^q(\Omega )} \le c(q,n,\Omega ) 
\Big(\|\bu^\ast\|_{\bW^{ -1,\, q}(\Omega )}   + \|g\|_{L^q(\Omega )}\Big).  
\ee
In particular, if $g=0$ we have $\bu\in W^{ 1,\, q}_{0, \rm div} (\Omega )^n$.
\end{thm}

Theorem\,\ref{thm5.7} shows the existence and uniqueness of weak $q$-solutions to Stokes system 
\eqref{5.1}, \eqref{5.2} holds for every bounded $ C^1$ domain.   Thus,   
according to Remark\,\ref{rem5.4}/2  we get the following 

\begin{cor}
\label{cor5.8}
Every bounded $C^1$-domain is $\nabla _q$-regular.  
\end{cor}

\begin{rem}
\label{rem5.9}
1. Let $1<q<+\infty $. The statement in Theorem\,\ref{thm5.7} continues to hold if $\Omega $ 
is a Lipschitz domain  with sufficiently small Lipschitz constant (cf. \cite{GaSiSo1994}).  In this case 
$\Omega $ is $\nabla _q$-regular. 

\hspace{0.5cm}
2.  Let $ \Omega \subset  \R^{3}$ be a bounded Lipschitz domain. As it has been proved  in \cite{BrSh1995}  there exists a 
number $3< q_0<+\infty$ such that the existence and uniqueness of weak $ q$-solutions 
holds for every  $q_0'< q<q_0$. Thus, $ \Omega $ is 
 $\nabla _q$-regular for all $q_0' < q < q_0$.  

\end{rem}

\hspace*{0.5cm} 
As a consequence of Theorem\,\ref{thm5.7} we get the existence and uniqueness of $q$-weak 
solutions to \eqref{5.1}, \eqref{5.2} for the case  $\delta =1$ too.

\begin{cor}
\label{cor5.10}
Let $\Omega \subset \R^n$ be a bounded $C^1$-domain. Let $1<q<+\infty $. 
Then for every $\bu^\ast\in  
W^{ -1,\, q}(\Omega )^n$    there exists a unique $q$-weak solution 
$\bu \in W^{ 1,\, q}_{0, \rm div} (\Omega )^n$ to \eqref{5.1}, \eqref{5.2} with $\delta =1$ and 
a unique 
$p\in  L^q_0(\Omega )$ satisfying 
\be\label{5.9}
-\Delta_q \bu + \bu= \bu^\ast-\nabla _q p \quad \mbox{in}\quad W^{ -1,\, q}(\Omega )^n.
\ee
In addition, there holds
\be\label{5.10}
\| \bu\|_{W^{ 1,\, q}(\Omega )}+ \|p\|_{L^q(\Omega )} \le c(q,n,\Omega) 
\|\bu^\ast\|_{\bW^{ -1,\, q}(\Omega )} .  
\ee
  
\end{cor}

{\bf Proof}:  In view of Remark\,\ref{rem5.2}/2.  we may restrict ourselves to the case $2\le q < +\infty $.  
In this case,  since every $q$-weak solution is a $2$-weak solution to \eqref{5.1}, \eqref{5.2}  the uniqueness is obvious.  

\hspace*{0.5cm} 
Let $\bu^\ast\in  W^{ -1,\, q}(\Omega )^n$  be arbitrarily chosen. 
Clearly, as $\bu^\ast\in  W^{ -1,\, 2}(\Omega )^n$ there exists a  $2$-weak solution to 
\eqref{5.1}, \eqref{5.2}. This solution  also is  a  $2$-weak solution to the Stokes system with right-hand side $\bu^\ast-\bu$.  Consulting Theorem\,\ref{thm5.7}  
replacing $\bu^\ast$ by $\bu^\ast-\bu$ therein  and applying 
a bootstrapping argument we obtain $\bu\in  W^{ 1,\, q}_{0, \rm div} (\Omega )^n$.  Whence, 
$\bu$ is a $q$-weak solution to  \eqref{5.1}, \eqref{5.2}. The estimate \eqref{5.10} follows from the closed mapping theorem.   
\hfill \Beweisende

\begin{rem}
\label{rem5.11}
Let $1<q<+\infty $.  If $\Omega $ is a bounded Lipschitz domain, 
then for every $f\in L^{q'} _0(\Omega )$ 
the equation $\nabla \cdot \bu=f$ has a solution in $W^{ 1,\, q'}_0(\Omega )^n$ (cf. 
\cite[Thm.III.3.1]{Galdi1994}). 
Hence, by the closed range theorem we see  that ${\rm im}\,\nabla _q= 
(W^{ 1,\, q}_{0, \rm div} (\Omega )^n)^\circ= G^{-1,\,q} (\Omega )^n$.  If   
$\Omega $ is  $\nabla _q$-regular,  then for every  $\bu^\ast\in  W^{ -1,\, q}(\Omega )^n$ 
we  may  define the  {\it associate pressure} as   
$p := \mathscr{P} (\bE^\ast_q \bu^\ast)$, which belongs to $L^q_0(\Omega )$. In addition, we have the estimate 
\be\label{5.11}
\|p\|_{L^q} \le c \|\bu^\ast\|_{W^{ -1,\, q}}, 
\ee  
where the constant $c>0$ depends on $q,n$, and the geometric properties of $\Omega $ only.  
\end{rem}

\hspace*{0.5cm} 
At the end of this subsection  we examine the scaling properties of the projection $\bE_q$. 
In particular, we will see that if $\Omega =B_R(x_0)$, the constant $c$ in \eqref{5.11} is  independent  of $R>0$.  
Without loss of generality we may assume that $x_0=0$, since all estimates which will be used are invariant under translations.  
Instead of $B_R(0)$ we shortly write  $B_R$.  Define, 
\[
\bPhi_q(\bu)(y) = R^{n/q'-1} \bu(Ry),\quad y \in B_1,\quad \bu \in W^{ 1,\, q'}_0(B_R)^n.   
\]
Clearly, according to  $\|\bPhi_q(\bu)\|_{W^{ 1,\, q'}_0(B_1)} =\|\bu\|_{W^{ 1,\, q'}(B_R)}$, the map  
$\bPhi_{q'} $ defines  an isometry between  $W^{ 1,\, q'}(B_R)^n$ and $W^{ 1,\, q'}(B_1)^n$. 
Accordingly,  its dual $\Phi ^\ast_q$ is an isometry between $ W^{ -1,\, q}(B_1)^n $ and $ W^{ -1,\, q}(B_R)^n$.  
On the other hand, it is readily seen that 
\[
\bPhi_{q'}  (W^{ 1,\, q'}_{0, \rm div} (B_R)^n) = W^{ 1,\, q'}_{0, \rm div} (B_1)^n,\quad 
\bPhi_{q}^\ast(G^{- 1,\,q}(B_1)^n )= G^{-1,\,q}(B_R)^n.  
\]
It can be easily checked,
\be\label{5.12}
\bE_{q', B_1} = \bPhi_{q'} \bE_{q', B_R} \bPhi_{q'}^{-1} \,\footnotemark\,.        
\ee 
\footnotetext{\,Obviously, the operator on the right is a projection with kernel $ W^{ 1,\, q'}_{0, \rm div} (B_1)^n $. For $\bfpsi\in C^\infty _{c}(B_1)^n$ we have $\bPhi_{q'}^{-1} (\bfpsi ) \in  C^\infty _{c}(B_1)^n$ and thus
$\bE_{q', B_1} (\bfpsi ) =\bPhi_{2} \bE_{2, B_R} \bPhi_{2}^{-1} (\bfpsi)$.  Since $\Phi _2$ is an isometry 
the operator on the right is self-adjoined with range $( W^{ 1,\, q'}_{0, \rm div} (B_1)^n )^{\bot} $. This shows that 
$\bPhi_{2} \bE_{2, B_R} \bPhi_{2}^{-1} (\bfpsi)= \bE_{2, B_1} (\bfpsi )$ and thus \eqref{4.6}. }
From \eqref{5.12} we infer 
\be\label{5.13}
\bE_{q, B_R}^\ast =   \bPhi_{q}^\ast \bE_{q, B_1}^\ast (\bPhi_{q}^\ast)^{-1}.        
\ee
In order to understand the relation between $\mathscr{P}_{B_R} $ and $\mathscr{P}_{B_1} $, we define 
\[
\Psi_q(f)(y) = R^{n/s} f(Ry),\quad y \in B_1,\quad f\in L^q(B_R).    
\]
Clearly, $\Psi _q$ defines an isometry between $L^q(B_R)$ and $L^q(B_1)$. 
Furthermore, we have $\Psi_q (L^q_0(B_R))= L^q_0(B_1)$. 
Let $p\in L^q(B_1)$.  By an elementary calculus we obtain 
\begin{align*}
\langle \bPhi^\ast_q \nabla_q p,  \bu \rangle &= - \intl_{B_1} p \nabla \cdot \bPhi_{q'} \bu   dy   
\\     
&= - R^{-n/s} \intl_{B_R} p(x/R) \nabla \cdot  \bu(x ) dx   = \langle \nabla_q \Psi ^{-1}_q p, \bu\rangle
\end{align*}
for all $\bu \in W^{ 1,\, q'}_0(B_R)^n$. Consequently, 
\be\label{5.14}
\bPhi _q^\ast \nabla _q = \nabla _q \Psi _q^{-1},\quad 
\mathscr{P}_{B_R}  \bE_{q, B_R}^\ast = 
\Psi _q^{-1} \mathscr{P}_{B_1} \bE_{q, B_1}^\ast (\Phi_q^\ast)^{-1}.  
\ee
Thus, we have the following 
\begin{cor}
\label{cor5.12}
There exists a constant $c=c(q,n)$ such that for all $0<R<+\infty $ the following is true: 
\be\label{5.15}
\|p\|_{L^q(B_R)} \le c \|\bu^\ast \|_{W^{ -1,\, q}(B_R)}\quad\forall\,   
\bu^\ast \in  W^{ -1,\, q}(B_R)^n,
\ee
where $p =  \mathscr{P}_{B_R}(\bE^\ast_{q, B_R}\bu^\ast) $.  
\end{cor}

\subsection{Exterior domains}

Next, let us investigate the  case  $\Omega $ being an exterior domain.  
For the existence and uniqueness 
of weak solutions to the Stokes-like system we argue similar as in \cite[Chap.\,V.]{Galdi1994}   with  slight modification. Accordingly, we have the following

\begin{thm}
\label{thm5.13}
Let $\Omega $ be an exterior $C^1$-domain.  Let $1<q<+\infty $. Then for every 
$\bu^\ast\in  W^{ -1,\, q}(\Omega )^n$  there exists a unique $q$-weak solution to 
\eqref{5.1}, \eqref{5.2} $(\delta =1)$.  In addition, there holds 
\be\label{5.16}
\|\bu\|_{W^{ 1,\, q}}\le c \|\bu^\ast\|_{W^{ -1,\, q}}.  
\ee
\end{thm}

\hspace*{0.5cm} 
Now, Theorem\,\ref{thm5.13} yields 

\begin{cor}
\label{cor5.14}
Every exterior $C^1$-domain $\Omega $  is $\nabla _q$-regular for all $1<q<+\infty $.  
\end{cor}

{\bf Proof}:  By virtue of Theorem\,5.13 we only need to show, that condition (b) in Lemma\,\ref{lem5.3} 
is fulfilled. Then the assertion will immediately follow from Lemma\,\ref{lem5.3}.  

{\it Proof of (b)}:  Let $\bfpsi \in C^\infty _{\rm c}(\Omega )^n $.  Then $\bu^\ast = 
-\Delta \bfpsi+ \bfpsi \in  C^\infty _{\rm c}(\Omega )^n $, which is embedded into 
$W^{ -1,\, q}(\Omega )^n$. According to Theorem\,\ref{thm5.13} there is a unique solution 
$\bu\in W^{ 1,\, q}_{0, \rm div} (\Omega )^n$ of \eqref{5.1}, \eqref{5.2} with $\delta =1$.  
By a standard  regularity argument making use of Theorem\,\ref{thm5.5} and  Corollary\,\ref{cor5.10}  we deduce that 
$\bu \in W^{ 1,\, 2}_{0, \rm div} (\Omega )^n$.  Whence (b), and the assertion of the corollary 
is completely proved. \hfill \Beweisende

%%%%%%%%%%%% 
% 
%section 6 
% 
%%%%%%%%%%%% 
\section{The associate pressure for time dependent distritions} 
\setcounter{secnum}{\value {section}
\setcounter{equation}{0}
\renewcommand{\theequation}{\mbox{\arabic{secnum}.\arabic{equation}}}}

Here we consider time distributions of the form $ \bu' + \bF$ 
in an interval $]a,b[$ 
$(-\infty <a<b<+\infty )$ (for the distributional time derivative cf. appendix below)  .   Let $\Omega \subset \R^n$ be a  domain. By $Q$ we denote the cylindrical domain $\Omega \times ]a,b[$.  

\hspace*{0.5cm} 
In the discussion below by  $L^s(a,b; L^q_{\rm loc} (\Omega ))$   we mean the linear space of all $f: ]a,b[ \rightarrow L^q_{\rm loc} (\Omega )$ such that 
\[
f|_{G\times ]a,b[}\in  L^s(a,b; L^q(G))\quad \forall\, G\Subset \Omega\quad  (1 \le s,q <+\infty).    
\]

\vspace{0.3cm}
Our first main result is the following

\begin{thm}
\label{thm6.1}
Let $1<q<+\infty $.  Let $\Omega \subset \R^n$ be a $\nabla _q$-regular domain together with the projection $\bE^\ast_q: W^{ -1,\, q}(\Omega )^n \rightarrow G^{ -1,\, q}(\Omega )^n $. 
Let $U \subset \Omega $ be a $q$-suitable subdomain.

\hspace*{0.5cm} 
1. Let $\bg, \bbf\in  L^1(a,b; W^{ -1,\,q}(\Omega )^n)$  such that 
\be\label{6.1}
\intl_{a}^{b}  \langle \bg(t), \bfpsi   \rangle \eta (t)  - \langle \bbf(t),\bfpsi  \rangle \eta'(t)dt  =0
\quad \forall\, \bfpsi \in  C_{\rm c, div}^\infty (\Omega )^n,\,\, \eta \in C^\infty _{\rm c}(]a,b[).
\ee
Then,  
\be\label{6.2}
\begin{cases}
{\D \intl_{a}^{b} \langle \bg(t), \bfvarphi(t) \rangle \eta (t)  - \Big \langle \bbf(t), \frac {d \bfvarphi } {d t}(t) \Big\rangle \eta'(t) dt  = 
\intl_{Q} (-p_{h} \partial _t \nabla \cdot \bfvarphi + p_{0}\nabla \cdot \bfvarphi  )   dx dt } 
\\[0.2cm]
\forall\, \bfvarphi \in  C_{\rm c}^\infty (Q)^n,
\end{cases}
\ee
where $p_{h}(t) = -\mathscr{P}^{(U)}_{ \Omega } (\bbf(t))$ and 
$p_{0}(t) =  -\mathscr{P}^{(U)}_{ \Omega }  (\bg(t))$  for a.\,e.  $t\in  ]a,b[$ (cf. Remark\,\ref{rem4.6}).

\vspace{0.2cm}
2.  If $\nabla \cdot \bbf =0$ in the sense of distribution, then $p_h(t) $ is harmonic for a.\,e. $t\in ]a,b[$.
If, in addition $\bbf \in  C_{w^\ast}([a,b]; W^{ -1,\, q}(\Omega )^n) $, then $p_h$ is continuous 
in $\Omega \times [a,b]$. 
\end{thm}
 
{\bf Proof}:  1. Since $C^\infty _{\rm  c, div} (\Omega )^n$ is dense in $W^{ 1,\, q'}_{0, \rm div} (\Omega )^n$ 
\eqref{6.1} remains true for all $\bfpsi \in  W^{ 1,\, q'}_{0, \rm div} (\Omega )^n = {\rm ker}\, \bE_{q'} $ 
(cf. \eqref{4.5} ).  Thus, we are in a position to apply Theorem\,\ref{thmA.4}  
for $X= W^{ 1,\, q'}_0(\Omega )^n, Y =X^{ \ast}= W^{ -1,\, q}_0(\Omega )^n $, 
$E= \bE_{q'} $, $f=\bbf, g^\ast=\bg$.  This, implies that 
$\bbf - \bE^\ast_q\bbf$ admits a distributive time derivative 
$(\bbf - \bE^\ast_q\bbf) ' = \bE_{q}^\ast \bg-\bg$ which satisfies  \eqref{A.12}, i.\,e.  
\be\label{6.3}
\intl_{a}^{b} \langle \bg(t), \bfvarphi(t)  \rangle - \Big\langle \bbf(t), \frac {d \bfvarphi } {d t}(t)
\Big\rangle dxdt  = 
\intl_{a}^{b} -\Big \langle \bE^\ast_q (\bbf(t)),  \frac {d \bfvarphi } {d t}(t) \Big \rangle  + \langle \bE_q^\ast (\bg(t)), \bfvarphi(t)  \rangle d t 
\ee
for all $\bfvarphi \in  C^\infty _{\rm c}(Q)^n $.  Define, 
\[
p_h(t) = - \mathscr {P}^{(U)}_{ \Omega } (\bg(t)),\quad  
p_0(t) = - \mathscr {P}^{(U)}_{ \Omega } (\bbf(t)),\quad 
\mbox{for a.\,e.\,\,$t\in ]a,b[$}. 
\] 
Clearly,   both  $ p_h$ and $ p_0$  belong to $L^1(a,b; L^q_{\rm loc}(\Omega))$ \,\footnotemark\,.  
Thus, \eqref{6.2} follows from \eqref{6.3} by using \eqref{4.15}. \hfill \Beweisende 
\footnotetext{\, Let $ \bbf \in L^1(a,b; W^{-1,\, q} (\Omega )^n)$. 
Observing \eqref{3.6} we see that for every ball $ B$, 
$\mathscr{P}^{(U)}_\Omega  (\bbf) |_{B} \in L^1(a,b; L^q(B))$, 
which shows that 
$\mathscr{P}^{(U)}_\Omega (\bbf)\in L^1(a,b; L^q_{\rm loc}(\Omega ))$. }

\vspace{0.2cm}
\hspace*{0.5cm}      
2.  If $\nabla \cdot \bbf=0$ in  the sense   of distributions,  then Lemma\,\ref{lem4.4} implies that  $p_h(t) $ is harmonic for a.\,e. $t\in  ]a,b[$.  

\hspace*{0.5cm} 
In addition, if $\bbf \in C_{w^\ast}([a,b]; W^{ -1,\, q}(\Omega )^n) $, from the first part of the theorem we see that 
\[
(\bbf+ \nabla_q  p_{h})' = (\bg  + \nabla_q p_{0}) \quad \mbox{in} \quad L^1(a,b; W^{-1,\, q}(G)^n).  
\]
Thus,  appealing to   Lemma\,\ref{lemA.1} with $X=Y= W^{ -1,\, q}(\Omega )^n $,  
we get $\bbf+ \nabla_q  p_{h} \in  C([a, b]; W^{ -1,\, q}(\Omega )^n)$. 
In particular, as $\bbf\in  C_{w^\ast} ([a,b]; W^{ -1,\, q}(\Omega )^n)$  we see that for every 
$\bfzeta \in  C^\infty_{\rm c}(\Omega )^n $  the function 
\[
t \mapsto \intl_{\Omega } \nabla  p_{h} (t) \cdot \bfzeta   dx = 
\intl_{\Omega } \langle \bbf(t),  \bE_{q'} \bfzeta \rangle  dx
\]
is  continuous on $[a, b]$. Taking $\bfzeta $ to be radial symmetric,  recalling  the mean value formula of harmonic functions it follows that
\begin{equation}
\nabla p_{h} (x_0, t) \rightarrow \nabla  p(x_0, t_0) \quad  \text{as} \quad  
t \rightarrow t_0 \quad \forall\,(x_0, t_0) \in  \Omega \times  [a,b].
\label{6.3a}
\end{equation}
 
\hspace*{0.5cm} 
Let $ (x_0,t_0) \in Q$ be fixed. Let $ B_R=B_R(x_0) \subset \Omega $.  
Recall that  $\mathscr{P}^{(U)}_\Omega |_{B_R}$ is a bounded linear operator from 
$W^{ -1,\, q}(\Omega)^n$ into $L^q(B_R)$ \,
(cf. Remark\,\ref{4.6} and \eqref{3.6}), which gives  
\[
\|p_{h}(t)\|_{L^q(B_R)}\le c \|\bbf(t)\|_{ W^{-1,\, q}(\Omega )}
\le c \| \bbf \|_{ L^\infty(a,b; W^{-1,\, q}(\Omega ))}\quad  \forall\,t \in [a,b].
\]
By using the properties of harmonic functions 
we find 
\[
\max_{x\in \overline{B_{R/2} }} |\nabla^2  p_{h}(x, t )|  \le  c R^{-n-2}  
\| \bbf \|_{ L^\infty(a,b; W^{-1,\, q}(\Omega ))}= C R^{ -n-2}\quad \forall\, t\in  [a,b]
\]
with a constant $C$ independent on $(x,t)$ and $ R$.  Using Newton-Leibniz formula we infer  
\[
|\nabla p_{h}(x, t ) - \nabla p_{h}(x_0, t )|\le CR^{ -n-2} |x-x_0|\quad 
\forall\, x\in B_{R/2}(x_0), \,\, t\in [a,b].
\]
Then by the aid triangle inequality we get 
\[
|\nabla p_{h}(x, t ) - \nabla p_{h}(x_0, t_0)|\le CR^{ -n-2} |x-x_0|+ 
|\nabla p_{h} (x_0, t )-\nabla p_{h}(x_0, t_0)|.
\]
Thus, taking into account \eqref{6.3a} we deduce  that  $\nabla p_{h}$ is continuous in 
$\Omega \times [a,b]$.    

\hspace{0.5cm}
As  $ p_h = p_h - (p_h)_U$ we are in a position to apply Lemma\,\ref{lemB.2}, 
which completes the proof of the theorem. 
\hfill \Beweisende

\begin{rem}
\label{rem6.5}
1. Let $\Omega \subset \R^n$ be any domain.  According to Definition\,\ref{def4.1}  $\Omega $ is 
$\nabla _2$-regular with the orthogonal projection $\bE^\ast_2$.  Thus, the statement of Theorem\,\ref{thm6.1} holds for $q=2$ without any restriction on  $\Omega $.  
\end{rem}

\vspace{0.3cm}
\hspace*{0.5cm} 
Next, we wish to introduce the local pressure projection associated to a bounded 
$C^1$-subdomain $G\subset \Omega $. To this end, we recall the definition of 
$\mathscr{P}_{q, G}: W^{ -1,\, q}(G)^n \rightarrow L^q_0(G)$, 
\[
\mathscr{P}_{q, G}(\bu^\ast) = \mathscr{P}_{q, G} ( \bE^\ast_{q, G} \bu^\ast),\quad 
\bu^\ast \in W^{ -1,\, q}(G)^n.
\]   
Let $1<q_1, q_2 <+\infty $. According to Remark\,\ref{rem3.4}, and Remark\,\ref{rem4.3}  we see that 
\be\label{6.4}
\mathscr{P}_{q_1, G}(\bu^\ast)=\mathscr{P}_{q_2, G}(\bu^\ast)\quad 
\forall\, \bu^\ast\in  W^{ -1,\, q_1}\cap W^{ -1,\, q_2} (G)^n. 
\ee
Hence, if no confusion can arise we omit the subscript $q$ and write 
$\mathscr{P}_{G} $ in place of $\mathscr{P}_{q, G}$.  Correspondently, we write 
$\bE^\ast_{G} $ in place of $\bE^\ast_{q, G}$.  Is readily seen that  
\be\label{6.5}
\nabla _q \mathscr{P}_{G}  = \bE^\ast_{G},\quad 
\mathscr{P}_{G} \nabla _q \mathscr{P}_{G}=\mathscr{P}_{G}. 
\ee
 
\vspace{0.3cm}
\hspace*{0.5cm} 
From consulting \cite{GaSiSo1994} we get  the following 

\begin{lem}
\label{lem6.2}
Let $G\subset \R^n$ be a bounded $C^{k} $-domain $(k\in \N)$. 
Let $1<q<+\infty $. Then $\mathscr{P}_{G} (\bbf) \in  W^{k-1,\, q}(G)^n $ for every 
$\bbf \in  W^{k-2,\, q}(G)^n$. In addition, there holds
\be\label{6.6}
\| \mathscr{P}_{G}(\bbf) \|_{W^{k- 1,\, q} (G)} \le c(q, k, n, G) 
\|\bbf\|_{\bW^{k-2,\, q}(G)}. 
\ee
Here $ W^{0,\, q}(G) $ equals $ L^q(G)$. 
\end{lem}

\vspace{0.3cm}
\hspace*{0.5cm} 
Our second main result  is the construction of the  following local pressure representation. 

\begin{thm}
\label{thm6.3}
Let $\Omega \subset \R^n$ be an open set.  Let $\bu\in  L^1(a,b; \bL^1_{\rm loc} (\Omega )) $   and $\bF\in  L^1(a,b; \bW^{ -1,\,q}_{\rm loc} (\Omega ))$  $(1<q<+\infty )$  such that 
\be\label{6.7}
\intl_{a}^{b} \langle \bF(t), \bfpsi   \rangle \eta (t)dt -
\intl_{Q} \bu(t)\cdot \bfpsi  \eta'(t)dxdt  =0
\quad \forall\, \bfpsi \in  C_{\rm c, div}^\infty (\Omega )^n,\,\, \eta \in C^\infty _{\rm c}(]a,b[).
\ee
Then, for every bounded subdomain $G\Subset \Omega $ with $C^1$-boundary we have 
\be\label{6.8}
\begin{cases}
{\D \intl_{a}^{b} \langle \bF(t), \bfvarphi(t)  \rangle dt - \intl_{Q} \bu \cdot \partial _t\bfvarphi dxdt  = 
\intl_{Q} (-p_{h, G} \partial _t \nabla \cdot \bfvarphi + p_{0,G}\nabla \cdot \bfvarphi  )   dx dt } 
\\[0.2cm]
\forall\, \bfvarphi \in  C_{\rm c}^\infty (G\times ]a,b[)^n,
\end{cases}
\ee
where $p_{h, G}(t) = -\mathscr{P}_{G} (\bu(t)|_{G} )$ and $p_{0,G}(t) =  
-\mathscr{P}_{G} (\bF(t)|_{G} )$  for a.\,e.  $t\in  ]a,b[$. 
\end{thm}

{\bf Proof}: Let $G\Subset \Omega $ be a bounded $C^1$-domain. 
Fix,  $1<s< \min \Big\{q, \frac {n} {n-1}\Big\}$. 
As $C^\infty _{\rm c, div}(G)^n $ is dense in $W^{ 1,\, s'}_{0, \rm div} (G)^n$,  \eqref{6.7} yields 
\be\label{6.9}
\intl_{a}^{b} \Big(\langle \bF(t), \bfpsi (t)  \rangle \eta(t) dt - 
\intl_{a}^{b} \intl_{G}  \bu(t)\cdot  \bfpsi dx\eta'(t)dt   =0
\ee
for all $\bfpsi \in W^{ 1,\, s'}_{0, \rm div} (G)^n= {\rm ker} (\bE_{s', G})  $ and $\eta \in C^\infty_{\rm c}(]a,b[)$.  
Thus, \eqref{6.9}  allows  to apply Theorem\,\ref{thmA.4} with $X= W^{ 1,\, s'}_{0} (G)^n$, $Y = L^1(G)^n$, and the projection $E = \bE_{s', G} $.  Indeed, 
by means of Sobolev's embedding theorem we have $Y \hookrightarrow X^\ast$  in the following sense 
\[
\langle \bv, \bfpsi  \rangle = \intl_{G} \bv\cdot \bfpsi  dx \le \|\bv\|_{L^1(G)}
\max_{\overline{G }}  |\bfpsi |  \le c\|\bv\|_{L^1(G)} \|\nabla \bfpsi \|_{L^{s'} (G)}
\]
 for all $\bv \in  L^1(G)^n$ and $\bfpsi \in  W^{ 1,\, s'}_{0} (G)^n $.    
Hence, the assumption \eqref{A.10}  of Theorem\,\ref{thmA.4} is fulfilled for  
$f(t) = \bu(t)$ and $g^\ast(t)= \bF(t)|_G $.  Consequently,  
\[
(\bu  - \bE_{s, G} ^\ast\bu)' + \bF  - \bE_{s, G} ^\ast \bF =0 \quad  
\mbox{in}\quad G\times ]a,b[
\]
in the sense of distributions.  Thus,  
setting $ p_{h, G}(t)= -  \mathscr{P}_{G} (\bu(t) |_{G}) $ and $p_{0, G}(t)= 
\mathscr{P}_{G}(\bF(t) |_{G}) $,   
from \eqref{A.12} we immediately get  
\[
\intl_{a}^{b} \Big\langle \bu(t)+ \nabla_q p_{h, G}(t), \frac {d \bfvarphi } {d t}(t)\Big\rangle   
 +   \intl_{a}^{b} \langle \bF(t) + \nabla_q  p_{0, G}(t), \bfvarphi (t) \rangle    d t=0 
\]
for all $\bfvarphi \in  C^1_{\rm c}(a,b; W^{1,\, s'}_{0}(G)^n)$. 
Whence, \eqref{6.8}. \hfill \Beweisende 

\vspace{0.3cm}
\hspace*{0.5cm} 
As consequence of Theorem\,\ref{thm6.3} we have 

\begin{cor}
\label{cor6.4}
Suppose all assumptions of Theorem\,\ref{thm6.3} are fulfilled.  Then the following statements are true. 

\vspace{0.3cm}
\hspace*{0.3cm} 1. Suppose $\nabla \cdot \bu =0$ in  the sense   of distributions. Then for every  bounded 
$C^1$-domain $G \Subset  \Omega $ the pressure $p_{h, G}(t) $ is harmonic for a.\,e. $t\in  ]a,b[$, and there holds 
\be\label{6.10}
\begin{cases}
{\D \intl_{a}^{b} \langle \bF(t), \bfvarphi(t)  \rangle dt - \intl_{Q} (\bu + \nabla p_{h, G}) \cdot \partial _t\bfvarphi dxdt  = \intl_{Q} p_{0,G}\nabla \cdot \bfvarphi   dx dt } 
\\[0.2cm]
\forall\, \bfvarphi \in  C_{\rm c}^\infty (G\times ]a,b[)^n.
\end{cases}
\ee
  
\hspace*{0.3cm}
2. If $\bu \in  L^1(a,b; \bL^q_{\rm loc} (\Omega ))$ for some $1<q<+\infty $ 
then for every   bounded 
$C^1$-domain $G \Subset   \Omega$  there holds $p_{h, G}\in  
L^1(a, b; W^{ 1,\, q}_{\rm loc} (G))$   fulfilling \eqref{6.8}.  If, in addition,  $ G$ is a $C^2$-domain 
we have  $p_{h, G}\in  L^1(a, b; W^{ 1,\, q} (G))$ together with the estimate 
\be\label{6.11}
\|\nabla p_{h, G}(t) \|_{\bL^q(G)} \le c(n,q, G) \|\bu(t)\|_{\bL^q(G)}\quad\mbox{for a.\,e. $t\in  ]a,b[$}.  
\ee

\hspace*{0.3cm}
3. 
If $\bF= \bF_1 + \ldots + \bF_N$, for $\bF_i \in L^1 (a,b; W^{ -1,\, q_i}_{\rm loc} (\Omega )^n )$ 
$(1< q_i < +\infty ; i=1,\ldots,N )$  then for every bounded $C^1 $-domain 
$G\Subset \Omega $ we have 
\be\label{6.12}
p_{0, G} = p_{0, G}^1 + \ldots + p_{0, G}^N, \quad \mbox{where}\quad p_{0, G}^i = 
\mathscr{P}_{G} (\bF_i|_{G}) \quad   
(i=1,\ldots,N).
\ee
Furthermore, there holds 
\be\label{6.13}
\| p^i_{h, G}(t) \|_{L^{q_i} (G)} \le c(n, q_i, G) \|\bF_i(t)\|_{W^{ -1,\, q_i}(G) } 
\quad\mbox{for a.\,e. $t\in  ]a,b[$} 
\ee
$(i=1,\ldots,N)$. 

\hspace*{0.3cm}
4. Assume that $\bu\in  C_{w} ([a,b]; L^1_{\rm loc} (\Omega )^n)$ with $\nabla \cdot \bu$ in  the sense   of distributions, then for every bounded $C^1$-domain $G\Subset \Omega $ 
the harmonic pressure $p_{h, G} $ is continuous in $G\times [a,b]$. 

\end{cor}

{\bf Proof}:  1.  As $\nabla \cdot \bu(t)=0 $   in  the sense   of distributions for a.\,e. $t\in  ]a,b[$ 
we may apply Lemma\,\ref{lem3.6} for $\bu(t)|_{G}$, which implies $p_{h, G } \in  
\overline{\mathscr{P}} (\bE^\ast_{q, G}  \bu(t)|_G) $ is harmonic for a.\,e. $t\in  ]a,b[$. 

\vspace{0.2cm}
\hspace*{0.3cm}
2. This first statement follows immediately from the local  regularity of weak solutions to the Stokes equation, 
while the second is an immediate consequence of \eqref{6.6} taking  $k=2$ therein 
(cf. Lemma\,\ref{lem6.2}). 

\vspace{0.2cm}
\hspace*{0.3cm}
3.  Set $s = \min \{q_1, \ldots, q_N\}$.  Owing to 
\[
p_{0, G}(t) = \mathscr{P}_{G} (\bF(t)|_{G})  =   \mathscr{P}_{G} (\bF_1(t) |_{G}) + \ldots + 
\mathscr{P}_{G} (\bF_N(t) |_{G}),
\]
the assertion is easily obtained by applying \eqref{6.6} (with $k=1$) to each 
of $p_{0, G}^i $ \, $(i=1,\ldots,N)$. 

\vspace{0.2cm}
\hspace*{0.3cm}
4. Let $G\Subset   \Omega $ be any bounded $C^1$-domain. According to the first statement of the corollary  
$p_{h, G} (t)$ is harmonic for every $t\in  [a,b]$, and by virtue of  Theorem\,\ref{thm6.3} we have 
\[
(\bu+ \nabla_q  p_{h, G})' = - \bF|_{G}  -\nabla_q p_{0, G}=0 \quad \mbox{in} \quad L^1(a,b; W^{-1,\, q}(G)^n).  
\]
Appealing to Theorem\,\ref{thm6.1}, we immediately see that 
$ p_{ h, G}= p_{ h, G}-(p_{ h, G})_G$ is continuous in $ G\times [a,b]$.
\hfill \Beweisende

%%%%%%%%%%%% 
% 
%section 7 
% 
%%%%%%%%%%%% 
\section{An application to distributional solutions to the generalized Navier-Stokes equation} 
\setcounter{secnum}{\value {section}
\setcounter{equation}{0}
\renewcommand{\theequation}{\mbox{\arabic{secnum}.\arabic{equation}}}}

We consider distributional solutions of the following generalized Navier-Stokes equations 
\begin{align}
\nabla \cdot \bu &=0 \quad \mbox{in}\quad Q,
\label{8.1}
\\
\partial _t \bu + \bu\cdot \nabla \bu - \nabla \cdot (a(x,t) \bD(\bu)) &= 
\nabla \cdot \bbf -\nabla p\quad \mbox{in}\quad Q,
\label{8.2}
\end{align}
where $\bu= (\bu^1, \ldots, \bu^n)$ denotes the velocity field ($n=2$ or $n=3$), 
$p$ the pressure, $a >0$ the viscosity  and $\nabla \cdot \bbf$ the external force. Here 
$\bD(\bu)$ stands for the matrix of the symmetric gradient given by 
\[
D_{ij}(\bu) = \frac {1} {2} (\partial _i u^j - \partial _j u^i),\quad i, j=1,\ldots,n. 
\]

Regarding  distributional solution to \eqref{8.1}, \eqref{8.2} we give the following definition
%\end{document}
\begin{defin}
\label{def7.1}
Let  $a\in  L^\infty (Q)$ and $\bbf\in L^2(Q)^n$. Then $\bu$ is called a 
{\it distributional solution with bounded energy} if 
\begin{align}
&\bu \in C_w([a, b]; L^2(\Omega )^n) \cap L^2(a,b; W^{ 1,\, 2}(\Omega )^n),\quad 
\label{8.3}
\\
&\nabla \cdot  \bu =0 \quad\mbox{a.\,e. in}\quad Q,
\label{8.4}
\end{align}
and the following identity holds for all $\bfvarphi \in  C^\infty _{\rm c}(Q)^n $ 
with $\nabla \cdot \bfvarphi =0$
\begin{align}
 &\intl_{Q} \bu \cdot \partial _t \bfvarphi - \bu\otimes \bu : \nabla \bfvarphi  + 
a\bD(\bu): \bD(\bfvarphi ) dxdt   
\cr
&\qquad= \intl_{Q} \bbf: \nabla  \bfvarphi   dx dt.
\label{8.5}
\end{align}

\end{defin}

\hspace{0.5cm}
The following result  is  a direct application of Theorem\,\ref{thm6.1},  and the results of Section\,6.

\begin{thm}
Given  $a\in  L^\infty (Q)$ and $\bbf\in L^2(Q)^n$, let $\bu$ be a distributional solution to \eqref{8.1}, \eqref{8.2} with bounded energy. Let $U\Subset \Omega $ be $2$-suitable. 
Then there holds 
\begin{align}\
 &\intl_{Q} (\bu +\nabla p_h)\cdot \partial _t \bfvarphi - \bu\otimes \bu : \nabla \bfvarphi  + 
a\bD(\bu): \bD(\bfvarphi ) dxdt   
\cr
&\qquad= \intl_{Q} \bbf: \nabla  \bfvarphi   dx dt + \intl_{Q} p_{0} \nabla \cdot \bfvarphi    dx dt 
\label{8.6}
\end{align}
for all $\bfvarphi \in C^\infty_{\rm c}  (Q)^n$, where 
\begin{align*}
 p_{h}(t)&= - \mathscr{P}^{(U)} ( \bE^\ast_2 \bu(t))  
\\     
p_{0}(t)&= 
 \mathscr{P}^{(U)} \Big(\bE^\ast (\nabla \cdot  (a\bD(\bu(t)) - \bu(t)\otimes \bu(t) + \bbf(t)))\Big)
\end{align*}
for a.\,e. $t\in  ]a,b[$. In particular, $p_h$ is harmonic with respect to $x\in  \Omega $ and continuous 
in $\Omega \times [a,b]$. 

\hspace*{0.5cm} 
If, in addition, $\Omega $ is a bounded $C^1$ domain then 
$p_0 = p_{0 }^1+ p_{0 }^2+ p_{0 }^3$, where 
\begin{align*}
p_{0}^1(t)&= 
 \mathscr{P} (\nabla \cdot  (a\bD(\bu(t))), 
\\     
p_{0}^2(t)&= -
 \mathscr{P}(\nabla \cdot (\bu(t)\otimes \bu(t))), 
\\
p_{0}^3(t)&= \mathscr{P}(\nabla \cdot \bbf(t)),
\end{align*}
while the harmonic pressure is given by $p_{h}(t)= -\mathscr{P}(\bu(t))  $ 
$(t\in  [a, b])$. 
\end{thm} 

{\bf Proof}:  Let $\bu $ be a distributional solution to \eqref{8.1}, \eqref{8.2}. 
 From the above definition it follows that $\bu \in  L^\infty (a,b; L^2(\Omega )^n)$.  
By a standard interpolation argument along with Sobolev's embedding theorem we infer 
\[
\bu \in  L^{8/n} (a, b; \bL^{4}(\Omega ) ).
\]
Thus, $a\bD(\bu) - \bu\otimes \bu + \bbf \in  L^{4/n} (a, b; L^{2}(\Omega)^n )$.  
Define, $\bF(t) \in  W^{ -1,\, 2}(\Omega )^n$ by 
\[
\langle \bF(t), \bv \rangle = \intl_{\Omega} (a\bD(\bu(t)) - \bu(t)\otimes \bu(t) + \bbf(t)): \nabla \bv  dx,
\quad \bv \in  W_0^{ 1,\, 2}(\Omega)^n.  
\]
for a.\,e. $t\in  ]a,b[$. As one can easily check $ \bF\in  
L^{4/n} (a, b; \bW^{ -1,\, 2}(\Omega)^n)$.  Now, applying Theorem\,\ref{thm6.1} together with Corollary\,\ref{cor6.4},  
  from \eqref{8.5} we get \eqref{8.6}.   
The second statement immediately follows from Corollary\,\ref{cor6.4}. 
\hfill \Beweisende

\appendix

%%% ----------------------------------------------------------------------
%       SECTION A
%%% ----------------------------------------------------------------------
\section{Vector valued functions}
\label{sec:-A}
\setcounter{secnum}{\value{section} \setcounter{equation}{0}
\renewcommand{\theequation}{\mbox{A.\arabic{equation}}}}

\hspace*{0.5cm}
Let $X$ be  a Banach space with norm  $\|\cdot \|_X$.  
Let $-\infty \le a < b \le +\infty $. 
By $L^s(a,b; X)$ ($1\le s \le +\infty$) we denote the   space of 
all Bochner measurable functions  $f: ]a,b[ \rightarrow X$ such that 
\[
\intl_{a}^{b} \|f(t)\|^s d t < \infty \,\,\, \mbox{if} \,\,1\le s< \infty; \quad
\esssup_{t\in (a,b)}  \|f(t)\|_{X} < \infty \,\,\,\mbox{if}\,\,s=\infty.   
\]

{\bf The Steklov mean}  Let $f \in  L^1(a,b; X)$. We extend $f$ outside $]a,b[$ by zero, and denote this extension again by $f$. For $\lambda \in  \R\setminus \{0\}$ we define the Steklov mean $f_\lambda: [a, b] \rightarrow X$, by 
\[
f_\lambda (t) = \frac {1} {\lambda } \intl_{t}^{t+\lambda } f(\tau ) d \tau,\quad t\in  [a,b] \,\footnotemark\,.   
\]
\footnotetext{\,Here, for $\lambda <0$ we set $\intl_{t}^{t+\lambda } f(\tau ) d \tau := 
-\intl_{t+\lambda }^{t} f(\tau ) d \tau $.}

The following properties of the Steklov mean are well-known and can be found in the standard literature:

\begin{itemize}
\item[(i)]  {\it $f_\lambda \in  C([a,b];X)$ for all $\lambda \in \R\setminus \{0\}$; }

\item[(ii)]  {\it If $f\in  L^s(a,b; X)$\,$(1\le s <+\infty )$ then}
\end{itemize}

\vspace*{-0.7cm}
\be\label{A.1}
\|f_\lambda \|_{ L^s(a,b; X)} \le  \|f \|_{ L^s(a,b; X)}\quad \forall\, \lambda \in \R\setminus \{0\},
\ee

\begin{itemize}
\item[]
{\it and $f_\lambda  \rightarrow f$ in $ L^s(a,b; X)$ as $\lambda  \rightarrow 0$ };

\item[(iii)] {\it Let $f\in  L^s(a,b; X)$ and $ g\in  L^{s'} (a,b; X)$\,$(1\le s\le+\infty )$, then}
\end{itemize}

\vspace*{-0.7cm}
\be\label{A.2}
\intl_{a}^{b} f_\lambda g  d t = \intl_{a}^{b} f  g_{-\lambda }   d t\quad \forall\, \lambda \in  \R\setminus 0.
\ee

Let $Y$ be a further  Banach space, its norm being denoted by $\|\cdot\|_Y$. 
Let $T \in \mathscr{L}(X,Y)$, i.\,e. $T: X \rightarrow Y$ is linear and bounded.  Then $T$  forms  a linear and bounded operator from 
$L^s(a,b; X)$  into  $ L^s(a,b; Y)$\,$(1\le s\le \infty)$, which again will be   
denoted by $T$ such that 
\be\label{A.3}
Tf(t) = T(f(t))  \quad \mbox{for a.\,e. $t\in ]a,b[$}. 
\ee
In particular,  we have 
\[
\intl_{a}^{b} \|Tf(t)\|_{Y}^s  d t\le  \|T\|^s  \intl_{a}^{b} \|f(t)\|_{X}^s  d t,\quad \forall\, f \in  L^s(a,b; X).
\]

{\bf Distributive time derivative}  Let $X, Y$ are Banach spaces such that $X$ is continuously, and densely embedded into $Y$.  Let $f\in  L^1(a,b; X)$. A Bochner function $g \in  L^1(a,b; Y) $  is called a {\it distributive time derivative }of $f$ if 
\be\label{A.4}
\intl_{a}^{b} g(t) \eta(t) d t = -\intl_{a}^{b} f(t) \eta'(t) d t \quad \mbox{in}\quad Y\quad \forall\, \eta \in  
C^\infty_{\rm c}(]a,b[).  
\ee 

\hspace{0.5cm}
Clearly, the distributive time derivative is unique and will be denoted by $f'$.  
The following important properties are well known and can be found in the standard literature. 

\begin{lem}
\label{lemA.1}
Let $f\in  L^1(a,b; X)$ with distributive time derivative $f' \in  L^1(a,b; Y)$. Then, 
eventually redefining $f$ on a subset of $[a,b]$ of measure zero we have 
$f\in  C([a,b]; Y)$, and there holds   
\be\label{A.5}
f(t) = f(s) + \intl_{s}^{t} f'(\tau )d \tau \quad \mbox{in}\quad  Y\quad \forall\, a\le s\le t\le b.  
\ee
In addition, we have 
\be\label{A.6}
\intl_{a}^{b} \langle v^\ast(t), f'(t)   \rangle d t = -
\intl_{a}^{b} \Big\langle \frac {d v^\ast } {d t} (t), f(t)  \Big\rangle   d t
\quad  \forall\,v^\ast \in C^1_{\rm c} ([a,b]; Y^\ast).
\ee

\end{lem}

\begin{lem}
\label{lemA.2}
Assume $X$ to be reflexive. Let $f \in  L^\infty (a,b;X)$ with distributive time derivative 
$f'\in L^1(a,b; Y)$. Then $f\in  C_w([a,b]; X)$.
\end{lem}

{\bf Proof}:  Thanks to  Lemma\,\ref{lemA.1} there holds $f\in  C([a, b]; Y)$.  Let $t\in  ]a,b[$. Then  
$(f_\lambda(t)) $ is bounded in $X$. Since $X$ is reflexive there exists a sequence $\lambda _j \rightarrow 0$ 
and $ \xi \in  X$ such that $f_{\lambda _j} (t) \rightharpoonup \xi $ in $X$ as $j \rightarrow +\infty $.  
According to $Y^\ast \hookrightarrow X^\ast$ and $f_{\lambda _j} (t) \rightarrow f(t)$ in 
$Y$ as $j \rightarrow +\infty $ we get 
\[
\langle v^\ast, f_{\lambda _j}(t)  \rangle \rightarrow \langle v^\ast, \xi  \rangle= \langle v^\ast, f(t)  \rangle 
 \quad \forall\, v^\ast \in  Y^\ast \quad  \text{as}\quad  j \rightarrow +\infty. 
\]
Consequently, $\xi = f(t)$.  In particular, $f(t) \in X$ for all $t\in  [a,b]$, and by the lower semi continuity of the norm we 
have 
\be\label{A.8}
\|f(t)\|_X\le \|f\|_{L^\infty (a,b; X)}\quad \forall\, t\in  [a,b].  
\ee

\hspace*{0.5cm}
Next, let $t_k \rightarrow t$ in $[a,b]$. Since $(f(t_k))$ is bounded in $ X$ (cf. \eqref{A.8}), as $X$ is reflexive there exists a subsequence 
$(t_{k_j} )$ and $\xi \in  X$ such that $f(t_{k_j} ) \rightharpoonup \xi $ as $j \rightarrow+\infty  $.   
As above, $f\in  C([a, b]; Y)$ yields $\xi = f(t)$. Since the limit is unique we get the convergence of the whole sequence, which proves the lemma. \hfill \Beweisende

\begin{lem}
\label{lemA.3}
Let $X,Y, Z$ are Banach spaces. As above we assume  $X \hookrightarrow Y$, densely. Let $T: Y \rightarrow Z$ be linear and bounded.  If $f\in  L^1(a,b; X)$ with distributive time derivative 
$f'\in  L^1(a,b; Y)$, then $Tf \in  L^1(a,b; Z)$ admits a distributive time derivative 
$(Tf)'\in  L^1(a, b; Z)$ given by 
\be\label{A.9}
(Tf)'(t) = Tf'(t) \quad \mbox{in}\quad Z \quad \mbox{for a.\,e. $t\in  ]a,b[$}.
\ee
\end{lem}

{\bf Proof}:  Clearly, with help of  \eqref{A.3} we get $ Tf' \in  L^1(a,b; Z)$.  
Let $\eta \in  C^\infty _{\rm c}(]a,b[)$.  Recalling the definition of the distributive 
time derivative, we see that 
\[
\intl_{a}^{b} f'(t) \eta(t) d t = -\intl_{a}^{b} f(t) \eta'(t) d t \quad \mbox{in}\quad Y.
\]
Then applying the operator $T$ to both sides of the  above identity,  we obtain 
\[
\intl_{a}^{b} Tf'(t) \eta(t) d t = -\intl_{a}^{b} Tf(t) \eta'(t) d t \quad \mbox{in}\quad Z.
\] 
This, shows that $(Tf)' = Tf'$, which completes the proof of the assertion. \hfill \Beweisende

\vspace{0.3cm}
{\bf Projections} Let $X, Y$ are Banach spaces such that  $X \hookrightarrow Y \hookrightarrow X^\ast $,  
each of the embeddings  are dense \,\footnotemark\,.  
\footnotetext{\,In the literature  $\{X, Y, X^\ast \}$ usual  is called a {\it Gelfand triple}.}

\vspace{0.3cm}
\hspace{0.5cm}
The following  theorem has been used  in the proof of our main result 
Theorem\,\ref{thm6.1}. 

\begin{thm}
\label{thmA.4}
Let $E: X \rightarrow X$ be a projection together with the direct sum   $X= {\rm im}\, E\oplus {\rm ker}\,E$.  
Let $f\in  L^1(a,b; Y)$ and $g^\ast \in  L^1(a,b; X^\ast)$ such that 
\be\label{A.10}
- \intl_{a}^{b}  \langle f(t), \psi \rangle \eta'(t) d t =
\intl_{a}^{b} \langle g^\ast(t), \psi \rangle \eta (t) d  t 
\ee
for all $\psi  \in {\rm ker} \,E$ and $\eta \in  C^\infty _{\rm c}(]a,b[) $.  
Then, $f-E^\ast f$ admits a distributive 
time derivative $ (f-E^\ast f)' \in  L^1 (a, b ; X^\ast)$ such that 
\be\label{A.11}
(f-E^\ast f)'(t) = g^\ast(t)- E^\ast g^\ast(t) \quad \mbox{in}\quad X^\ast\quad \mbox{for a.\,e. $t\in ]a,b[$}\,\footnotemark\,. 
\ee
\footnotetext{\, Here $E^\ast : X^\ast \rightarrow X^\ast$ stands for  the dual of $E$. 
}
In particular, there holds 
\be\label{A.12}
- \intl_{a}^{b}  \Big\langle f(t)- E^\ast f(t), \frac {d \varphi } {d t}(t)\Big\rangle d t =
\intl_{a}^{b} \langle g^\ast -  E^\ast g^\ast(t), \varphi (t) \rangle d  t 
\ee
for all $\varphi \in C^1_{\rm c} (a,b; X)$.

\end{thm}

\vspace{0.3cm}
{\bf Proof}:  From the assumption of the theorem it follows that  
$ I- E$ is  a projection from $X$ onto ${\rm ker}\, E$.  From 
\eqref{A.10} we deduce 
\[
- \intl_{a}^{b}  \langle f(t), (I- E)\psi \rangle \eta'(t) d t =
\intl_{a}^{b} \langle g^\ast(t), (I-E)\psi \rangle \eta (t) d  t 
\]
for all $\psi \in X$ and $\eta \in  C^\infty _{\rm c}(]a,b[) $.  This  shows that 
\[
- \intl_{a}^{b}  (f(t)-E^\ast f(t)) \eta'(t) d t =
\intl_{a}^{b} (I-E^\ast) g^\ast(t) \eta (t) d  t \quad  \text{in}\quad  X^{ \ast}\quad \forall\, \eta \in  C^\infty _{\rm c}(]a,b[).  
\]
Thus, $(f -E^\ast f)' \in  L^1(a,b; X^\ast)$, and there holds  \eqref{A.11}. 
Finally, \eqref{A.12} can be obtained 
with help of  \eqref{A.6}  taking into account the canonical embedding $X \hookrightarrow X^{\ast \ast} $. \hfill \Beweisende

%%%%%%%%%%%% 
% 
%section B 
% 
%%%%%%%%%%%% 
\section{Continuity of potentials of continuous  gradient fields} 
\setcounter{secnum}{\value {section}
\setcounter{equation}{0}
\renewcommand{\theequation}{\mbox{B.\arabic{equation}}}}

In this appendix we discuss the question of continuity of time dependent gradient fields being continuous in a subcyliner.  This result has been used  
in the proof of the continuity of the harmonic pressure (cf. Theorem\,\ref{thm6.1}, Corollary\,\ref{cor6.4}). 

\begin{lem}
\label{lemB.1}
Let $\Omega \subset  \R^n$ be a domain. Let $-\infty <a<b<+\infty $.  
Let $ p : \Omega \times [a,b] \rightarrow \R$ such that 
$ p(\cdot, t)\in C^1(\Omega )$ for all $t\in  [a,b]$. 
Furthermore, suppose that 
$\nabla p $  is continuous in $\Omega \times [a,b]$, and  
there exists $x_0\in \Omega $ such that $p(x_0, \cdot )\in C([a,b])$. Then 
$ p $  is continuous in $\Omega \times [a,b]$. 
\end{lem}

{\bf Proof}:  Firstly, we claim that $p(x, \cdot ) \in C([a,b])$ for all $x\in \Omega $. 
To prove this, let us denote by $\Omega _c$  the set of all points $x\in \Omega  $ such that 
$p(x, \cdot ) \in C([a,b])$.  Since $x_0 \in \Omega _c$, this set is nonempty. Thus, in order 
to prove that $\Omega _c=\Omega $ we only need to show that $\Omega _c$ is open and relatively 
closed. 

\vspace{0.3cm}
{\it (i) $\Omega _c$ is open}. Let $x \in  \Omega_c $. Fix a ball $B\Subset  \Omega $  centered  in $x$.  Then  by applying the Newton-Leibniz formula 
for all $y\in  B$ we find
\begin{align*}
p (y, t)&= p(y, t) - p(x,t) +p (x, t)   
\\     
&=  \intl_{0}^{1} \nabla p(x + \tau (y-x),t )\cdot (y-x) d \tau  
+p(x, t).
\end{align*}
Due to our assumptions both functions on the right-hand side belong to $C([a,b])$, and thus $x\in  \Omega _c$. 

\vspace{0.3cm}
{\it (ii) $\Omega _c$ is relatively closed} Let $x \in \overline{\Omega}^{\rm rel} _c \subset \Omega $. 
Let  $B\Subset  \Omega $ be a ball having its center  in $x$. Clearly,  there exists 
$y\in B\cap  \Omega_c $. As above we see that 
\begin{align*}
p(x, t)&= p(x, t) - p (y,t) +p(y, t)   
\\     
&=  \intl_{0}^{1} \nabla p(y + \tau (x-y),t )\cdot (x-y) d \tau  
+p(y, t).
\end{align*}
Hence, as the term on the right-hand side belongs to $C([a,b])$ we  deduce that $x\in \Omega _c$.  Whence, $\Omega _c =\Omega $. 

\vspace{0.3cm}
\hspace*{0.5cm} 
Secondly, let $(x,t)\in  \Omega \times [a,b]$. Let $B\Subset \Omega $ be a ball with center $x$. 
By using the triangle inequality, and Newton-Leibniz formula, we get for all $(y,s)\in B\times [a,b] $
\[
|p(x,t)-p(y,s)|\le |p(x,t)-p(x,s)|+ \intl_{0}^{t} \nabla p (y+\tau (x-y), s) \cdot (x-y) d \tau.   
\]
From this inequality  we infer that $p$ is continuous in $(x,t)$, since the first term tends to zero 
as $s \rightarrow t$ according to the first part of the proof, while the second term tends to $ 0$ as $ y \rightarrow x$, 
since $\nabla p(\cdot , s)$ is bounded on $B$. \hfill \Beweisende

\begin{lem}
\label{lemB.2}
Let $\Omega \subset  \R^n$ be a domain. Let $-\infty <a<b<+\infty $.  
Let $ p : \Omega \times [a,b] \rightarrow \R$ such that 
$ p(\cdot, t)\in C^1(\Omega )$ for all $t\in  [a,b]$. Furthermore, suppose that 
$\nabla p $  is continuous in $\Omega \times [a,b]$. Then for every subdomain 
$U\subset \Omega $ with $\mes U<+\infty $  and $(p- p_{U})  |_{U}\in L^\infty(a,b; L^q(U)) $ \,
$(1<q<+\infty )$ the function  $p - p_U $ is continuous in $\Omega \times [a,b]$. 
\end{lem}

{\bf Proof}:  According to Lemma\,\ref{lemB.1} it will be sufficient to prove the existence of $x_0\in \Omega $ 
such that $p(x_0, \cdot )- p_U\in C([a,b])$.  If $U\Subset \Omega $ is a ball with center $x_0$, then 
by  using Newton-Leibniz formula it follows
\[
p (x_0,\cdot )- p(\cdot )_U = \frac {1} {\mes U} \intl_{U} \intl_{0}^{1} \nabla p (x-\tau (x_0-x), \cdot)\cdot (x-x_0) d \tau   dx.
\]
Since $\nabla p $ is continuous on $\overline{U }\times [a,b]$ we infer that 
$p(x_0, \cdot ) - p(\cdot )_U\in C([a,b])$. Consequently, $p-p_{U} $ is continuous in $\Omega \times [a,b]$. 

\hspace{0.5cm}
Next, suppose $U\Subset \Omega $.  Let $B\Subset \Omega $ be a ball with center $x_0$. Then 
\[
p(x_0,\cdot )- p(\cdot )_{U} =  (p(x_0, \cdot)- p(\cdot)_{B})+ (p(\cdot)_{B}- p(\cdot)_{U}). 
\]
The function which occurs in first parenthesis belongs to $C([a,b])$, which  has been shown above,  while  the second one can be evaluated as follows 
\[
p(\cdot)_{B}- p(\cdot)_{U} = \frac {1} {\mes U}\intl_{U} p(\cdot )_{B}- p(y, \cdot )  dy.   
\]
As the integrant belongs to $C(\overline{U }\times [a,b])$ the function on the left-hand side is continuous on $[a,b]$. 
Thus, Lemma\,\ref{lemB.1} implies the assertion. 

\vspace{0.3cm}
\hspace{0.5cm}
Finally, if $U\subset \Omega $ with $\mes U <+\infty $ and $ (p- p_{U})  |_{U}\in L^\infty(a,b; L^q(U)) $ 
we may choose a sequence $U_m \Subset U_{m+1}\Subset  U$ of  increasing domains such that 
$\cup _{m=1}^\infty U_m =U $.  According to the second step we have 
$p - p_{U_m} = \pi - \pi _{U_m} $ is continuous in $\Omega \times [a,b]$, 
where $\pi = p-p_{U} $. 
Let $x_0\in  \Omega $.  As $\pi _{U}=0 $ in order to prove that $p(x_0, \cdot ) - p(\cdot) _{U} = 
\pi (x_0,\cdot ) $ is continuous on $[a,b]$ it will be sufficient to verify that $\pi _{U_m} \rightarrow 0 $ uniformly on $[a,b]$. 
In fact, this is true since 
\[
\pi(t) _{U_m}= - \frac {1} {\mes U_m} \intl_{U\setminus U_m} \pi  (y, t)dy \le 
\frac {\mes(U\setminus U_m)^{1/q'} } {\mes U_1}\|\pi \|_{L^\infty (a,b; L^q(U))}.   
\]
Hence, Lemma\,\ref{lemB.1} gives the desired continuity of $p-p_{U} $. \hfill \Beweisende

\vspace{0.5cm}
{\bf Acknowledgements}
The author has been supported by the Brain Pool Project of the Korea Federation of Science and Technology Societies  (141S-1-3-0022).  
The author   also wishes to thank Prof. H.-O. Bae for his meaningful  comments and suggestions  on the present paper, as well as  
the hospitality during the visit  in Korea at Ajou university, Suwon.

%\bibliographystyle{siam} 
%\bibliography{LocPressProj}

\end{document}